%% file: main.tex
\documentclass[10pt,journal,compsoc]{IEEEtran}
\IEEEoverridecommandlockouts
\usepackage[noadjust]{cite}
\usepackage{cite}
\usepackage{amsmath,amssymb,amsfonts,dsfont}
\usepackage{algorithm} 
\usepackage{algpseudocode} 
\usepackage{graphicx}
\usepackage{textcomp}
\usepackage{xcolor}
\usepackage{comment}
\usepackage{adjustbox}
\usepackage{float}
\usepackage{hyperref}

\makeatletter
\long\def\@makecaption#1#2{\ifx\@captype\@IEEEtablestring%
\footnotesize\begin{center}{\normalfont\footnotesize #1}\\
{\normalfont\footnotesize\scshape #2}\end{center}%
\@IEEEtablecaptionsepspace
\else
\@IEEEfigurecaptionsepspace
\setbox\@tempboxa\hbox{\normalfont\footnotesize {#1.}~~ #2}%
\ifdim \wd\@tempboxa >\hsize%
\setbox\@tempboxa\hbox{\normalfont\footnotesize {#1.}~~ }%
\parbox[t]{\hsize}{\normalfont\footnotesize \noindent\unhbox\@tempboxa#2}%
\else
\hbox to\hsize{\normalfont\footnotesize\hfil\box\@tempboxa\hfil}\fi\fi}
\makeatother

\usepackage{dblfloatfix}
\usepackage[caption=false,labelformat=empty]{subfig}
\DeclareMathOperator*{\argmax}{arg\,max}
\DeclareMathOperator*{\argmin}{arg\,min}

\newcommand{\E}{\mathbb{E}}
\newcommand{\R}{\mathbb{R}}
\newcommand{\N}{\mathbb{N}}
\newcommand{\Z}{\mathbb{Z}}
\newcommand{\Char}{\mathds{1}}
\newcommand{\ones}{\mathbf{1}}

\newcommand{\vx}{\mathbf{x}}
\newcommand{\vy}{\mathbf{y}}

\newcommand{\cS}{\mathcal{S}}

\newcommand{\cN}{\mathcal{N}}

\newcommand{\floor}[1]{\!\left\lfloor #1 \right\rfloor\!}
\newcommand{\ceil}[1]{\!\left\lceil #1 \right\rceil\!}
\newcommand{\pibf}{\boldsymbol{\pi}}

\newtheorem{theorem}{Theorem}[section]

\newtheorem{corollary}[theorem]{Corollary}
\newtheorem{proposition}[theorem]{Proposition}
\newtheorem{lemma}[theorem]{Lemma}

\newtheorem{remark}[theorem]{Remark}
\newcommand{\qedsymbol}{\hfill\square}

\begin{document}

\title{Minimizing File Transfer Time in Opportunistic Spectrum Access Model 
%-----AUTHOR INFO AND THANKS-------
\author{Jie~Hu,
        Vishwaraj~Doshi,
        and~Do~Young~Eun% <-this % stops a space
\IEEEcompsocitemizethanks{\IEEEcompsocthanksitem 
A subset of the material in this paper was presented in \cite{hu2021opportunistic}.\\
Jie Hu and Do Young Eun are with the Department of Electrical and Computer Engineering, and Vishwaraj Doshi is with the Operations Research Graduate Program, North Carolina State University, Raleigh, NC. Email: \{jhu29, vdoshi, dyeun\}@ncsu.edu. This work was supported in part by National Science Foundation under Grant CNS-1824518 and IIS-1910749.\protect\\}
}

}

\IEEEtitleabstractindextext{%
%----------ABSTRACT------------ 
\begin{abstract}
We study the file transfer problem in opportunistic spectrum access (OSA) model, which has been widely studied in throughput-oriented applications for max-throughput strategies and in delay-related works that commonly assume identical channel rates and fixed file sizes. Our work explicitly considers minimizing the file transfer time for a given file in a set of heterogeneous-rate Bernoulli channels, showing that max-throughput policy doesn't minimize file transfer time in general. We formulate a mathematical framework for static extend to dynamic policies by mapping our file transfer problem to a stochastic shortest path problem. We analyze the performance of our proposed static and dynamic optimal policies over the max-throughput policy. We propose a mixed-integer programming formulation as an efficient alternative way to obtain the dynamic optimal policy and show a huge reduction in computation time. Then, we propose a heuristic policy that takes into account the performance-complexity tradeoff and consider the online implementation with unknown channel parameters. Furthermore, we present numerical simulations to support our analytical results and discuss the effect of switching delay on different policies. Finally, we extend the file transfer problem to Markovian channels and demonstrate the impact of the correlation of each channel.
\end{abstract}

\begin{IEEEkeywords}
Opportunistic spectrum access, file transfer problem, minimum transfer time, shortest path problem
\end{IEEEkeywords}
}

\maketitle

\IEEEdisplaynontitleabstractindextext

% \begin{center} Cover letter \end{center}
% \fi

\IEEEpeerreviewmaketitle

%---------SECTIONS---------
\input{Sections/Introduction}
\input{Sections/Problem_description}
\input{Sections/Analysis_static_case}
\input{Sections/MDP_formulation}

\input{Sections/Practical_concern}
\input{Sections/Simulation}

\bibliographystyle{IEEEtran}
\bibliography{ref}

\vspace{-15mm}
\begin{IEEEbiography}[{\includegraphics[width=1in,height=1.25in,clip,keepaspectratio]{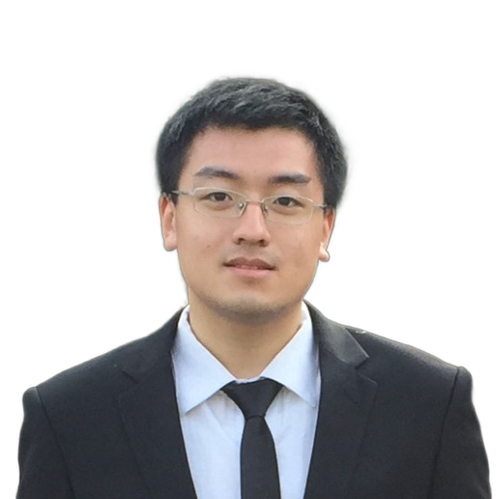}}]%
{Jie Hu}
received his B.E degree in communication engineering from Wuhan University of Technology, Wuhan, China, and Masters degree in electrical engineering from Northwestern University, Evanston, IL, USA. He is a Ph.D. student in the Department of Electrical and Computer Engineering at North Carolina State University. His current research interests are in the area of machine learning in dynamic spectrum access problem.
\end{IEEEbiography}
\vspace{-15mm}
\begin{IEEEbiography}[{\includegraphics[width=1in,height=1.25in,clip,keepaspectratio]{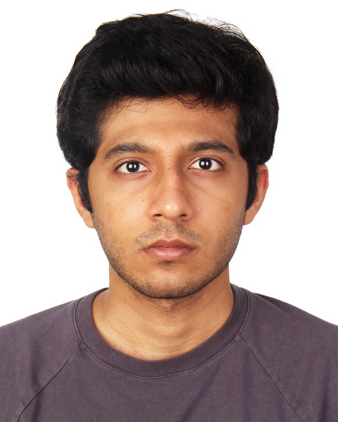}}]%
{Vishwaraj Doshi}
received his B.E. degree in mechanical engineering from the University of Mumbai, Mumbai, MH, India, and Masters degree in Operations Research from North Carolina State University, Raleigh, NC, USA. He is currently pursuing his Ph.D. degree with the Operations Research Graduate Program at North Carolina State University. His primary research interests include design of randomized algorithms on graphs, and epidemic models on networks.
\end{IEEEbiography}
\vspace{-10mm}
\begin{IEEEbiography}[{\includegraphics[width=1in,height=1.25in,clip,keepaspectratio]{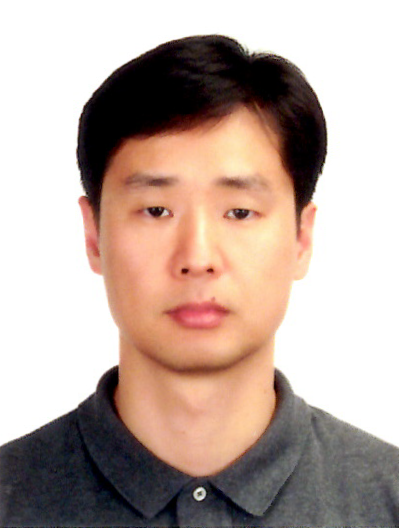}}]%
{Do Young Eun}
(Senior Member, IEEE) received his B.S. and M.S. degree in Electrical Engineering from Korea Advanced Institute of Science and Technology (KAIST), Taejon, Korea, in 1995 and 1997, respectively, and Ph.D. degree from Purdue University, West Lafayette, IN, in 2003. Since August 2003, he has been with the Department of Electrical and Computer Engineering at North Carolina State University, Raleigh, NC, where he is currently a professor. His research interests include distributed optimization for machine learning, machine learning algorithms for networks, distributed and randomized algorithms for large social networks and wireless networks, epidemic modeling and analysis, graph analytics and mining techniques with network applications. He has been a member of Technical Program Committee of various conferences including IEEE INFOCOM, ICC, Globecom, ACM MobiHoc, and ACM Sigmetrics. He is serving on the editorial board of IEEE Transactions on Network Science and Engineering, and previously served for IEEE/ACM Transactions on Networking and Computer Communications Journal, and was TPC co-chair of WASA'11. He received the Best Paper Awards in the IEEE ICCCN 2005, IEEE IPCCC 2006, and IEEE NetSciCom 2015, and the National Science Foundation CAREER Award 2006. He supervised and co-authored a paper that received the Best Student Paper Award in ACM MobiCom 2007. 
\end{IEEEbiography}

\end{document}

%% file: Sections/Introduction.tex
\section{Introduction}\label{Introduction}

In recent years, there has been an explosion in demand for wireless services due to the rapid growth in the number of wireless devices, including mobile devices and Internet of Things (IoT) devices. This demand further exacerbates the scarcity of allocated spectrum, which is ironically known to be underutilized by licensed users \cite{osa_survey}. The \textit{Opportunistic spectrum access} (OSA) model has been proposed to reuse the licensed spectrum in an opportunistic way otherwise wasted by licensed users \cite{osa_survey}. Recently, the FCC has released a new guidance in $2020$, which would expand the ability of the unlicensed devices (especially IoT devices) to operate in the TV-broadcast bands \cite{federal_2020}. Besides, the related IEEE 802.22 family has been developed to enable spectrum sharing \cite{802.22_standard} to bring broadband access to rural areas.

In the OSA model, a \textit{secondary user} (SU) aims to opportunistically access the spectrum when it is not used by any other users, while also prioritizing the needs of the \textit{primary user} (PU). The SUs need to periodically sense the spectrum to avoid interfering with PUs. We call an SU's behavior \textit{static} if it adheres to only one channel, and \textit{dynamic} if it is free to switch channels. While the concept of this model is simple, the design of spectrum sensing strategy faces various challenges: the interaction among multiple secondary users \cite{naparstek2018deep,avner2019multi,nasir2019multi}, spectrum sensing policy in the Markovian channels \cite{liu2010indexability,yadav2022deep}, the trade-off between the cost of sequential sensing (when permitted) and the expected reward \cite{Online_Sequential_2,zuo2019obp}, channel selection under resource constraints \cite{monemian2015optimum,gan2018cost,liang2019deep}, to list a few.

\subsection{Motivation: Throughput v.s. Latency}\label{motivation}
Nowadays, low latency has become one of the main goals for $5$G wireless networks \cite{Low_Latency_Towards_5G} and other time-sensitive applications with guaranteed delay constraints. In many applications, data is valid only for a limited duration and should be delivered before it expires, and vehicular communication is one such scenario. The increased demand for intelligent vehicular traffic (e.g., autonomous car development \cite{hussain2018autonomous}) has led to the need for vehicular communication to explore spectrum holes for offloading vehicular users in device-to-device mode \cite{paul2019spectrum}. In addition, delay-sensitive safety messages (i.e., speed and position of the vehicles) require low latency (as low as 100 ms) \cite{8024608}. Another example is medical body sensor networks \cite{sodagari2018technologies}, where the cognitive radio is implemented in body sensor network for life-critical monitoring, e.g., packets indicating a patient's health abnormality should be sent to a doctor as soon as possible, especially when the patient is out of the hospital network and needs to temporarily borrow vacant spectrum resources.

In the OSA literature, \emph{throughput} is one of the most commonly used performance metrics. Recent studies \cite{Srikant_2018,Srikant_2019,zhu2020machine,almasri2021managing}, by utilizing the multi-armed-bandit (MAB) techniques, have focused on finding max-throughput channel while the SU needs to learn the unknown channel parameters on the fly. In order to transmit a file as quickly as possible, common folklore might assume that the max-throughput policy would also suggest the minimum expected file transfer time. For example, Wald's equation implies that the file download time in an \emph{i.i.d} (over time) channel is equal to the file size divided by the average throughput of that channel, implicitly favoring the max-throughput channel for minimal download time. However, as will be explained in Section \ref{static analysis}, we find that this is \emph{not} the case in general. 

On the other hand, most \textit{delay}-related works consider average queuing delay of a large number of packets (fixed size) following Poisson arrivals \cite{wu2014learning,cao2017dynamic, dimitriou2018stable, huang2019dynamic,iqbal2021enhanced}. However, focusing on individual file transfers is important for small files when the SU needs to transmit each file as soon as possible. For instance, IEEE 802.11p protocol requires each car to generate and send safety messages continuously at $100$ ms intervals \cite{8024608}, making Poisson arrivals unsuitable to model this situation. In addition, the file size can vary depending on the application \cite{tschabitscher_2020}. Same channel data rate across all channels is another implicit assumption in those delay-related works \cite{wu2014learning,cao2017dynamic, dimitriou2018stable, huang2019dynamic,iqbal2021enhanced},\footnote{Channel data rate differs from the service rate in \cite{wu2014learning,cao2017dynamic, dimitriou2018stable, huang2019dynamic} because service rate is related to the length of time the channel is available, while data rate refers to the speed of data transfer through a channel.} but it doesn't reflect the realistic heterogeneous channel environment assumed in the throughput-oriented studies \cite{Srikant_2018,Srikant_2019,zhu2020machine,almasri2021managing}. Clearly, allowing the SU to \emph{switch} over such channels during instances of PU’s interruption can further reduce the file transfer time, but to the best of our knowledge, this issue has not been fully explored.

\subsection{Related Works and Their Limitations}\label{status_quo}
\emph{Throughput} and \emph{delay} are two performance metrics commonly used in the OSA literature to evaluate the quality of service (QoS) in the wireless network. For throughput-oriented works, the PU's behavior can be modeled as a two-state Markov chain (thus correlated over time), for which partially observable Markov decision processes (POMDPs) are typically employed to formulate the spectrum sensing strategy in order to maximize the long-term throughput \cite{Zhao_2007}. These POMDPs do not possess known structured solutions in general and they are known to be Polynomial-Space-complete (PSPACE-complete) even if all the channel statistics are known a priori \cite{liu2015online}. To achieve near maximum throughput, computationally efficient yet sub-optimal policies, such as myopic policy \cite{zhao2008opportunistic} and Whittle's index policy \cite{liu2010indexability}, have been proposed for the \emph{offline} OSA setting (known channel parameters). In particular, both \cite{zhao2008opportunistic} and \cite{liu2010indexability} introduced a concept of ``belief vector'' to guess the available probability of each channel and updated the vector after each observation of the channel state. In each time slot, the myopic policy \cite{zhao2008opportunistic} selected the channel with the maximum ``guess'' throughput, while the Whittle's index policy \cite{liu2010indexability} selected the channel with the highest value according to the Whittle's index and the belief vector.

Recently, machine learning techniques have emerged in the \emph{online} setting (unknown channel parameters) that the SU needs to learn the unknown channel environment in order to find the max-throughput policy. For example, model-based MAB techniques \cite{tekin2011online,dai2012efficient,NEURIPS2019_2edfeadf} and model-free deep neural networks \cite{wang2018deep,liang2019deep}, are utilized to obtain the max-throughput policy over Markovian channels. Explicitly, single-channel online policies have been developed in \cite{tekin2011online,dai2012efficient} to find the channel with maximum long-term throughput. \cite{NEURIPS2019_2edfeadf} utilized thompson sampling to estimate the parameters of each channel and employed the \textit{offline} policies from \cite{zhao2008opportunistic,liu2010indexability}. With well-trained neural networks, \cite{wang2018deep} showed better performance than the policies in \cite{zhao2008opportunistic,liu2010indexability}. \cite{liang2019deep} tackled the coordination problem among multiple SUs with deep Q-learning. On the other hand, MAB techniques have been extensively studied for heterogeneous channels, each with \textit{i.i.d} Bernoulli distribution, in order to find the channel with maximum throughput. The widely used MAB techniques include the Bayesian approach \cite{mohamedou2017bayesian}, upper confidence bounds \cite{almasri2021managing}, thompson sampling \cite{Srikant_2018} and its improvement from efficient sampling \cite{Srikant_2019}, and coordination approach among multiple SUs \cite{naparstek2018deep,avner2019multi}. In both two channel models studied in the throughput-oriented works, i.e., Bernoulli channels and Markovian channels, we will later show that ``the max-throughput channel does not always minimize the file transfer time''.

For delay-sensitive applications, packet delay in cognitive radio networks has been extensively studied using queuing theory to derive delay-efficient spectrum scheduling strategies. In this setting, a stream of packet arrivals (of the same size) modeled as a Poisson process with a constant rate is a common assumption in delay related works \cite{wu2014learning,cao2017dynamic, dimitriou2018stable, huang2019dynamic,iqbal2021enhanced}, and the goal is often to minimize the average packet delay in the steady state. Specifically, a spectrum sensing strategy (including queuing delay, packet priority and interruption by PU's) was studied in the multimedia applications \cite{wu2014learning}.  A dynamic load-balancing spectrum decision scheme was proposed in \cite{cao2017dynamic}, where an R-learning algorithm was introduced to deal with unknown channels and queuing statistics. \cite{dimitriou2018stable} controlled the transmission probability of each SU in the random access network and proposed a random-access strategy for two-and-three-SUs cases to minimize the queuing delay. Later, a model-free reinforcement-learning based strategy was studied in \cite{huang2019dynamic} to predict the channel accessing schemes without information exchanges among SUs and maximize the Quality-of-Service performance of the target SU. \cite{iqbal2021enhanced} focused on the hybrid spectrum access strategy with interweave and underlay spectrum access techniques for multiple SUs in order to obtain both good throughput and lower packet delay. However, as discussed Section \ref{motivation}, fixed packet length with Poisson arrivals and the same channel rate assumed in \cite{wu2014learning,cao2017dynamic, dimitriou2018stable, huang2019dynamic,iqbal2021enhanced} is not realistic in some delay-sensitive applications. The file transfer problem considered in this paper allows for arbitrary file sizes and heterogeneous channel rates, and we can tackle the file transfer problem for each single file.

\subsection{Our Contributions}
In this paper, we study the OSA model with the aim of minimizing the transfer time of a single file over the heterogeneous Bernoulli channels with different channel rates and channel available probabilities. We also provide practical implementations that take into account computational costs and unknown channel environments, as well as the extension to Markovian channels. Our main contributions can be summarized as follows:

\smallskip
\noindent \textit{Theoretical Analysis of the File Transfer Problem.}
    \begin{itemize}
        \item We first analyze the expected file transfer time of a single file for \textit{static} policies.
        \item By using the analysis of the static policy as a stepping stone, we interpret the file transfer problem under \textit{dynamic} policies as a stochastic shortest path problem, and obtain the dynamic \textit{optimal} policy.
        \item We show that both static and dynamic optimal policies reduce the transfer time compared to the baseline max-throughput policy, and this reduction is even more significant in delay-sensitive applications, where files are relatively small.
    \end{itemize}
    
\noindent \textit{Practical Considerations.}
    \begin{itemize}
        \item By formulating a mixed-integer programming method, we present an alternative technique to solve the shortest path problem, which speeds up the computation of dynamic optimal policy.
        \item We propose a lightweight heuristic policy to further reduce the computational cost while maintaining good performance compared to the max-throughput policy.
        \item In the online setting, where channel parameters are unknown to the SU, we modify an MAB algorithm proposed in \cite{talebi2017stochastic} and show its gap-dependent regret bound, guaranteeing the learning of the optimal policy that minimizes the file transfer time.
    \end{itemize}
     
\noindent \textit{Simulations.} 
    \begin{itemize}
        \item  We empirically show that the max-throughput policy is not the best when it comes to achieving the minimum file transfer time in both known and unknown channel environments. 
        
        \item When channel switching delay is taken into account, our lightweight heuristic policy can even outperform the dynamic optimal policy.
    \end{itemize}

\noindent \textit{Extension to Markovian Channels.} 
    \begin{itemize}
        \item We extend the theoretical analysis to Markovian channels, where each channel is modeled as a two-state Markov chain, and obtain the expected file transfer time of a single file in each channel and show the effect of correlation on the file transfer time.
    \end{itemize}

The rest of the paper is organized as follows: In Section \ref{model description} we introduce the OSA model and characterize the file transfer problem and its policy under the OSA framework. In Section \ref{static analysis}, we show the expected transfer time for static policy and it's performance analysis. Then, we extend from the static policy to the dynamic policy in Section \ref{mdp formulation and dynamic policy}. The practical concerns are discussed in Section \ref{online learning}. In Section \ref{simulation}, we evaluate different policies in the numerical setting. We provide additional analysis on the file transfer problem over Markovian channels in Section \ref{File Transfer in Markovian Channels}.

%% file: Sections/Problem_description.tex
\section{Model description}\label{model description}
\subsection{The OSA Model}\label{problem setting}
Consider a set of $N$ heterogeneous channels $\cN \triangleq \{1,2,\cdots,N\}$ available for use and each channel $i\in\cN$ offers a stable rate of $r_i>0$ bits/s if successfully utilized \cite{dai2012efficient,liu2010indexability}. In our setting, a SU wishes to transfer a file of size $F$ bits using one of these $N$ channels via opportunistic spectrum access. The SU can only access one channel at any given time, and can maintain this access for a fixed duration of $\Delta$ seconds, after which it has to sense available channels again (even the same channel) in order to avoid the interference to the active PU (or other SUs) in the current channel. At this point, the SU can decide which channel to sense and access that channel for the next $\Delta$ second interval if the channel is \textit{available}. Or the SU has to wait for $\Delta$ seconds to sense again if that channel is \textit{unavailable}, thereby unable to transfer data for this duration. This pattern is known as the \textit{constant access time} model, and has been commonly adopted for the SU's behavior as a collision prevention mechanism in the OSA literature \cite{osa_survey,mohamedou2017bayesian}. The cycle repeats itself until the SU transmits the entire file size $F$, then it immediately exits the channel in use. We omit the channel switching delay in our OSA model and the duration $\Delta$ seconds are fully used for file transmission, which is typically assumed in order to simplify the mathematical model and design a throughput-optimal policy in the OSA literature \cite{liu2010indexability,Srikant_2018,Srikant_2019,zhu2020machine,almasri2021managing,tan2022cooperative}. Note that the duration $\Delta$ seconds is not a randomly chosen number. For example, $\Delta$ is recommended as $100$ ms because the SU needs to vacate the current channel within $100$ ms once the PU shows up, as defined in IEEE 802.22 standard \cite{802.22_standard}. The SU can transmit up to $3.1$ Mb in each $\Delta$ seconds with highest channel rate $31$ Mbps in IEEE 802.22 standard and many small files (e.g, $5$ KB text-only email, $800$ KB GIF image and $3$ MB YouTube short video) need just a few slots to transmit.

\begin{remark}
The duration $\Delta$ does not include the sensing time for the fair comparison between our policies in Section \ref{static analysis} to \ref{online learning} and the max-throughput policy in the same OSA model  \cite{Srikant_2018,Srikant_2019,zhu2020machine,almasri2021managing}. In addition, as simulated in \cite{pei2009much}, if the time duration is set to $100$ ms (which is the duration of our $\Delta$) and the target probability of accurate sensing is around $90\%$, the sensing time for a cognitive radio network is typically chosen to be $6$ ms. This sensing time is negligible compared to the whole time duration and is omitted in our mathematical model.
\end{remark}

We say a channel is unavailable (or \textit{busy}) if it is currently in use by the primary users (PUs) or other SUs, while it is available (or \textit{idle}) if it is not in use by any other users. The state of a channel (idle or busy) is assumed to be independent over all channels $i\in\cN$, and \emph{i.i.d.} over the time instants $\{0,\Delta,2\Delta,\cdots \}$ following Bernoulli distribution with parameter $p_i\in(0,1]$, in line with the widely used discrete-time channel model \cite{osa_survey,mohamedou2017bayesian,avner2019multi}.\footnote{We also extend our theoretical analysis of the file transfer problem over Markovian channels in Section \ref{File Transfer in Markovian Channels}, i.e., each channel is modeled as a two-state Markov chain that will change its state accordingly every $\Delta$ seconds.} Specifically, for each $i \in \cN$, $\{Y_i(k)\}_{k\in\N}$ is a Bernoulli process with $p_i=P[Y_i(k)=1]=1-P[Y_i(k)=0]$ for all $k\in\N$.
Then, we can define $X_i(t)$, the state of channel $i\in\cN$ at any time $t \in \R_+$, as a piecewise constant random process $X_i(t) \triangleq Y_i(\floor{\frac{t}{\Delta}})$,
where $\floor{\cdot}$ denotes the floor function. This way, we write $X_i(t) = 1$ (or $0$) if channel $i\in \cN$ is available (or unavailable) for the SU with probability $p_i$ (or $1-p_i$). 

\begin{remark}
Inaccurate sensing, including mis-detections with probability $v_i$ for channel $i$ in each time slot, has been studied in the OSA literature \cite{pei2009much,toma2020estimation}. However, this does not affect our theoretical analysis. With probability $p_i' \triangleq p_i (1- v_i)$, the SU can successfully transmit data in the current time slot, otherwise the data transmission is zero. Thus, we can take the inaccurate sensing into account and replace $p_i$ with $p_i'$ in the analysis in Section \ref{static analysis} to \ref{online learning} without affecting the conclusions.
\end{remark}

The rate at which the SU can transmit files through channel $i\in\cN$ at any time instant $t\geq0$, also termed as the \textit{instantaneous throughput} of the channel $i$, is given by $r_iX_i(t)$, with its \textit{throughput} \cite{Srikant_2019} by $\E[r_i X_i(t)] = r_ip_i$. We denote by $i^* \triangleq \argmax_{k\in\cN} r_kp_k$ the channel with the \textit{maximum throughput}. For simplicity, we assume that this channel is unique, i.e., $r_{i^*}p_{i^*} \!>\! r_k p_k$ for all $k \in \cN \!\setminus\! \{ i^*\}$. In the next section, we take a closer look at the max-throughput policy and static policies in general.

\subsection{Policies for File Transfer}\label{subsection:policies for file transfer}

\begin{figure}[ht]\vspace{-2mm}
    \centering
    \includegraphics[width=0.85\columnwidth]{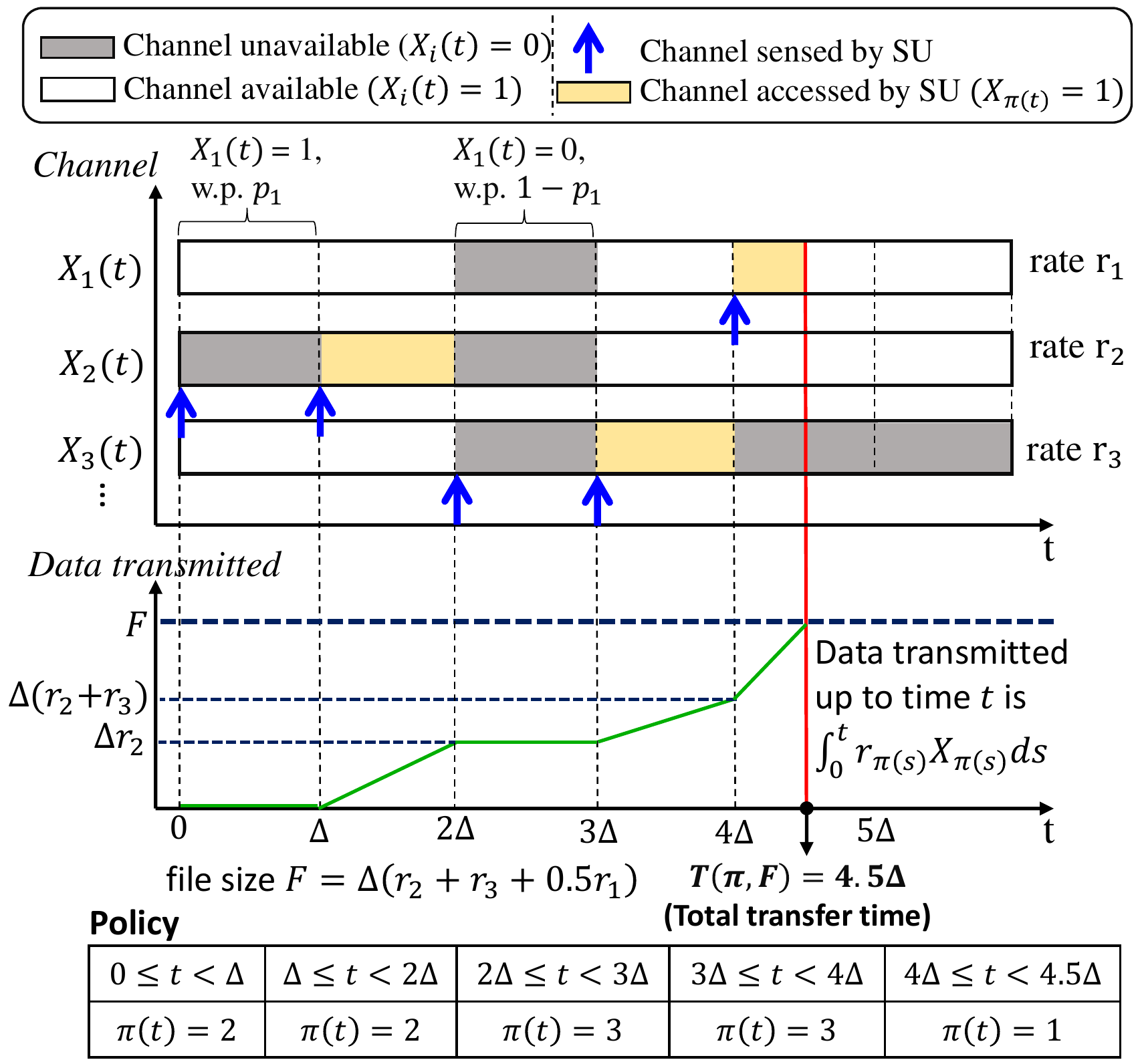} 
    \vspace{-2mm}
    \caption{File transfer via the OSA framework. The SU senses channels according to its policy, and accesses the channel if it is available. Upon gaining access, it begins transmitting data at the corresponding channel rate. Transmission ends (red line) as soon as the amount of data transmitted (green line) equals the file size $F$.}
    \label{fig:illustration}\vspace{-1mm}
\end{figure}

We define a \textit{policy} at time $t$ to be a mapping $\pi: \R_+ \mapsto \cN$ where $\pi(t)=i$ indicates that the SU has chosen channel $i$ to access during the time period $\left[\floor{\frac{t}{\Delta}}\!\Delta, \left(\floor{\frac{t}{\Delta}} \!\!+\!\!1 \right)\!\Delta\right)$. From our standing assumption, a policy therefore only changes at $t \!\in\! \{0,\Delta,2\Delta,\cdots\}$, and all policies ensure that the file transfer for any finite size $F$ will eventually be completed. This way, the policy $\pi(t)$ is a piecewise constant function (mapping), defined at all time $t \geq 0$. For a given policy $\pi$, let $T(\pi,F)$ denote the \emph{transfer time} of a file of size $F$ --- the entire duration of time to complete the file transfer, which is written as
%-----------------
\begin{equation}\label{stopping time}
    T(\pi,F) = \min_{T\geq 0}\left\{\int_0^{T} r_{\pi(t)} X_{\pi(t)}(t)dt ~~\geq F \right\}.
\end{equation}
Figure \ref{fig:illustration} explains the file transfer progress via OSA model.
%-----------------
%----------------------------

The objective of our OSA framework is to minimize the expected transfer time $\E[T(\pi,F)]$ over the set of all policies $\pi$. Policies can be \textit{static}, where the SU only senses and transmits via one (pre-determined) channel $i$ throughout the file transfer, that is, $\pi(t) = i$ for all $t\geq0$. For such static policies, we denote by $T(i,F)$ their transfer time for file size $F$. The channel that provides the minimum expected transfer time is then called \textit{static optimal} given by
%-----------------
\begin{equation}\label{eqn:static_optimal_channel}
    i_{so}(F) = \text{arg min}_{k \in \cN} ~\E [T(k,F)], \quad \text{(static optimal)}
\end{equation}
%-----------------
\noindent and we denote $T(i_{so},F)$ the corresponding transfer time for this static optimal policy. Note that the static optimal channel $i_{so}(F)$ depends on the file size $F$ and can vary for different file sizes. Policies can also be \textit{dynamic}, in which an SU is allowed to change the channels it chooses to sense throughout the course of the file transfer. Given a file size $F$, the policy with the
minimum expected transfer time over the set of all policies $\Pi(F)$ is called the \textit{dynamic optimal} policy given by \begin{equation}\label{eqn:dynamic_optimal_channel}
    \pi^*(F) \!=\! \text{arg min}_{\pi \in \Pi(F)} ~\E [T(\pi,F)].~ \text{(dynamic optimal)}
\end{equation} 
Lastly, we define the \textit{max-throughput policy} as the static policy with the channel $i^*$, which maximizes the long-term throughput. In the next section, we take a closer look at the max-throughput policy and static policies in general.

%% file: Sections/Analysis_static_case.tex
\section{Static optimal policy}\label{static analysis}

Recent works in the OSA literature focus on estimating channel parameters $p_i$'s, with the goal of eventually converging to the policy $i^* \!=\!\! \argmax_{k\in\cN} \! r_kp_k$ which provides the maximum throughput \cite{mohamedou2017bayesian,dai2012efficient,Srikant_2018,Srikant_2019,zhu2020machine,almasri2021managing}.\footnote{While \cite{Srikant_2018,Srikant_2019} deal with link rate selection problem to select best rate in one channel to maximize the expected throughput, the mathematical model of link rate selection problem is essentially the same as the standard OSA setting for choosing the max-throughput channel, as considered in our setting.} They focus on minimizing the `regret' in the MAB model, defined as the difference between the cumulative reward obtained by the online algorithm and the max-throughput policy (the optimal policy in hindsight).

The essential assumption behind all these approaches is that the SU always fully dedicates $\Delta$ seconds in each time interval for file transfer. Channel $i^*$ appears as a good candidate since it provides the largest expected data transfer $\Delta r_{i^*} p_{i^*}$ across every time interval. This is further supported by the well-known Wald's equation with the \textit{i.i.d} reward assumption at each time interval, suggesting that $\E[T(i,F)] \!=\! F/r_ip_i$ for each channel $i\!\in\!\cN$, which is then minimized by $i^*$. When policies are dynamic, however, the rewards are not identically distributed since the transfer rates of the dynamically accessed channels can be different, making Wald's equation inapplicable. Surprisingly, it is not applicable for static policies either. As typically is the case in delay-sensitive applications \cite{wu2014learning,dimitriou2018stable, huang2019dynamic,iqbal2021enhanced}, the file sizes are often not that large, rendering their transfer times small enough that an SU may not need to utilize the whole $\Delta$ seconds for data transfer in each time interval. The reward summands are still not identically distributed, causing Wald's equation to be inapplicable in general.

Our key observation in this paper is that choosing channel $i^*$ may not be the best option to minimize the expected transfer time. In this section, we limit ourselves to the set of static policies of the form shown in \eqref{eqn:static_optimal_channel} and analyze the resulting expected transfer time in the OSA network. We use this to compare the performance gap between the max-throughput policy and the static optimal policy, and show that for a reasonable choice of channel statistics and file sizes, the static optimal policy performs significantly better than the max-throughput policy. 
We derive a closed-form expression of the expected transfer time of a file of size $F$ in each fixed channel by the following proposition.

%------------------------------------------
%----------PROPOSITION---------------------
%------------------------------------------
\begin{proposition}\label{prop:expected_time_static} 
Given a file of size $F$, the expected transfer time $\E[T(i,F)]$ of the static policy for channel $i \in \cN$ is
\begin{equation}\label{eqn:time_static_policy}
     \E[T(i,F)] = \Delta \left(k_i/p_i + \Char_{\{\alpha_i > 0\}}(1-p_i)/p_i + \alpha_i\right),
\end{equation}
where $ k_i \triangleq \floor{F/\Delta r_i} \in \Z_+$ and $\alpha_i \triangleq  F/ \Delta r_i - k_i \in [0,1)$.
\end{proposition}
%------------------------------------------
%----------PROOF---------------------------
%------------------------------------------
\begin{IEEEproof}
Observe that any file size $F$ transmitted in channel $i$ can be written as
%------------------------------------------
\begin{equation}\label{eqn:file_size_static}
    F = k_i \Delta r_i + \alpha_i \Delta r_i
\end{equation}
%------------------------------------------
with $k_i$ being the number of intervals fully utilized for successful transmission, and $\alpha_i$ being the fraction of the $\Delta$ second interval utilized for file transfer toward the end. After choosing channel $i$, the SU first spends a random amount of time, denoted by $T_{wait}^n$ ($n = 1,2,\cdots,k_i$), waiting for channel $i$ to become available and starts the $n$-th transmission in that channel for $\Delta$ seconds. If the remaining portion $\Delta \alpha_i$ is not zero, the SU needs additional random waiting time $T_{wait}^{k_i+1}$ to complete the transfer. These random variables $\{T_{wait}^n\}$ are geometrically distributed and i.i.d over $n=1,2,\cdots,k_i+ \Char_{\{\alpha_i>0\}}$ with mean $E[T_{wait}^n] = \Delta (1-p_i)/p_i$. Let the constant $T_{tran} \triangleq F/r_i = \Delta(k_i+\alpha_i)$ be the total successful transmission time. Then the transfer time can be written as 
\begin{equation*}
    \begin{split}
        T(i,F) &= \sum_{n=1}^{k_i} \left(T_{wait}^n + \Delta\right) + \Char_{\{\alpha_i > 0\}} T_{wait}^{k_i+1} + \Delta \alpha_i \\
        &= T_{tran} + \!\!\sum_{n=1}^{k_i+ \Char_{\{\alpha_i > 0\}}}\!\!  T_{wait}^n.
    \end{split}
\end{equation*}
Taking the expectation of the equation above yields \eqref{eqn:time_static_policy}.
\end{IEEEproof}

\smallskip
\noindent
From Proposition \ref{prop:expected_time_static}, the expected transfer time $\E[T(i,F)]$ for any file size $F$ under the static policy on channel $i \in \cN$ can be explicitly written in terms of file size $F$, time duration $\Delta$ and channel statistics $r_i$ and $p_i$ of the chosen channel $i$. Substituting $k_i = F/\Delta r_i - \alpha_i$ in \eqref{eqn:time_static_policy} gives 
\begin{equation}\label{eqn:static_lower_bound}
    \E[T(i,F)] \!=\! \frac{F}{r_ip_i} \!+\! \Delta \Char_{\{\alpha_i>0\}}  (1\!-\!\alpha_i) \frac{1\!-\!p_i}{p_i} \!\geq\! \frac{F}{r_ip_i}.
\end{equation}
The inequality in \eqref{eqn:static_lower_bound} shows the expected transfer time of any static policy is no smaller than that given by Wald's equation.

\begin{figure}[ht]
    \centering
    \vspace{-0mm}
    \includegraphics[width=0.85\columnwidth]{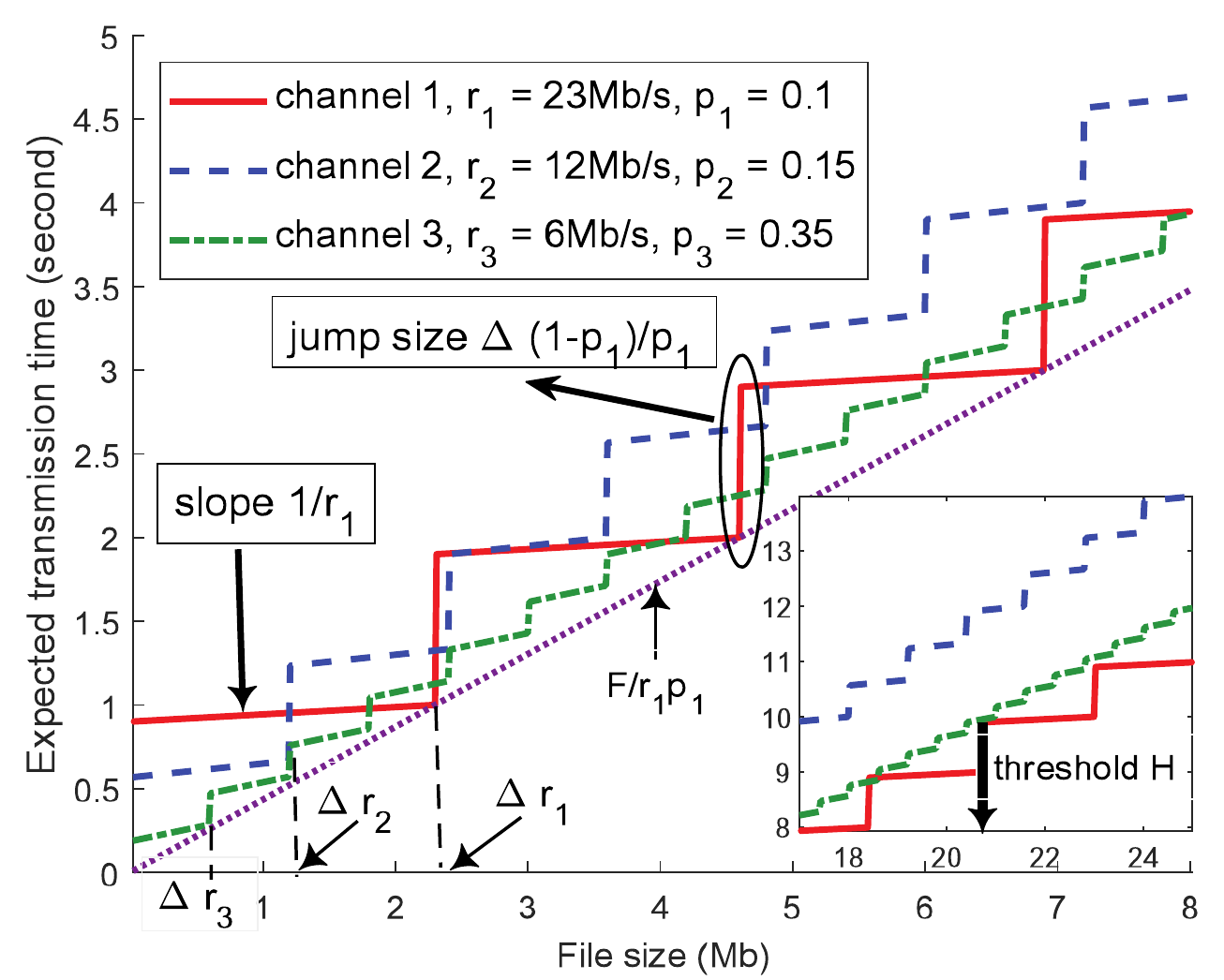}
    \vspace{-2mm}
    \caption{Expected transfer time from \eqref{eqn:time_static_policy} with duration $\Delta = 100$ ms (IEEE 802.22 standard). Channel $1$ has the maximum throughput. The purple dot-line $E[T] = F/r_1p_1$ corresponds to the Wald's equation for channel $1$. Downward arrow in the inset is the threshold $H$ in Proposition \ref{prop:threshold_optimal_policy}.}
    \label{fig:time_3channel_static}
    \vspace{-1mm}
\end{figure}

We use Figure \ref{fig:time_3channel_static} to illustrate the results in Proposition \ref{prop:expected_time_static}, where each line represents the expected transfer time via one channel over a range of file sizes from \eqref{eqn:time_static_policy}. We observe that the expected transfer time of channel $1$ (red line) is always above the Wald's equation of channel $1$ (purple dot-line). As shown in \eqref{eqn:time_static_policy}, the slope of each line (channel $i$) is $1/r_i$ and the ``jump size'' $\Delta(1-p_i)/p_i$ is equal to the expected waiting time till the channel is available. The jumps in the plot for each channel $i$, representing the waiting times, occur at exactly the instances where file size is an integer multiple of $\Delta r_i$, and come into play especially when there is still a small amount of remaining file to be transferred at the end of a $\Delta$ time interval.
%------------------------------------
%------------------------------------

By definition, the static optimal policy provides the minimum expected transfer time over all static policies including the max-throughput policy itself. While it is true for all file sizes, in some cases with certain file sizes, these two policies may coincide.

%------------------------------------------
%----------PROPOSITION---------------------
%------------------------------------------
\begin{proposition}
\label{prop:threshold_optimal_policy}
The max-throughput policy coincides with the static optimal policy, that is, $i_{so}(F) = i^*$, for any file size $F$ satisfying at least one of the two conditions below:
\begin{enumerate}
    \item $F$ exceeds a threshold $H$, where
%------------------------------------------
\begin{equation} \label{eqn:static_policy_threshold}
    H = \frac{\Delta (1-p_{i^*})/p_{i^*}}{1/r_hp_h - 1/r_{i^*}p_{i^*}},
\end{equation}
and $h = \argmax_{j \in \cN \setminus \{i^*\}} r_jp_j$ is the channel with the second largest throughput.
%------------------------------------------
    \item $F$ is an integer multiple of $\Delta r_i^*$, i.e., $F \!=\! k \Delta r_{i^*}$ for some $k \!\in\! \Z_+$. 
\end{enumerate}
\end{proposition}
%------------------------------------------
%----------PROOF---------------------------
%------------------------------------------
\begin{IEEEproof}
From $\Char_{\{\alpha_i>0\}}(1-\alpha_i) \in [0,1]$ and \eqref{eqn:static_lower_bound}, we have the upper bound and the lower bound of $E[T(i,F)]$ as follows:
% %------------------------------------------
\begin{equation}\label{eqn:lower_upper_bound}
\begin{aligned}
    \frac{F}{r_ip_i} \leq \E[T(i,F)] &= \frac{F}{r_ip_i} + \Delta \Char_{\{\alpha_i>0\}}  \frac{(1-\alpha_i) (1-p_i)}{p_i} \\ &\leq \frac{F}{r_ip_i} + \Delta \frac{1 - p_i}{p_i}.    
\end{aligned}
\end{equation}
To ensure $\E[T(i^*,F)]  \leq \E[T(j,F)]$ for all $j \in \cN \setminus \{ i^*\}$, it suffices to consider the upper bound of $E[T(i^*,F)]$ to be always smaller than the lower bound of $E[T(j,F)]$ for all $j \in \cN \setminus \{ i^*\}$ from \eqref{eqn:lower_upper_bound}, that is 
\begin{equation}\label{eqn:sufficient_inequality}
\begin{aligned}
\frac{F}{r_{i^*}p_{i^*}} + \Delta \frac{1-p_{i^*}}{p_{i^*}} \leq \min_{j \in \cN \setminus \{ i^*\}}\frac{F}{r_jp_j}.
\end{aligned}
\end{equation}
By definition of channel $h = \argmax_{j \in \cN \setminus \{i^*\}} r_jp_j$ we have $F/r_hp_h = \min_{j \in \cN \!\setminus\! \{ i^*\}}F/r_jp_j$. Then rearranging the second inequality in \eqref{eqn:sufficient_inequality} yields $F \geq H$ in (a).

When $F=k \Delta r_{i^*}$, we have $\alpha_{i^*} = 0$. Then from \eqref{eqn:time_static_policy}, the expected transfer time is simply $\E[T(i^*,F)] = \Delta\left(k/p_{i^*}\right) = F/r_{i^*}p_{i^*}$. Since $r_{i^*}p_{i^*} \geq r_jp_j$ for any $j\in\cN$, we have
$$\E[T(i^*,F)] = \frac{F}{r_{i^*}p_{i^*}} \leq \frac{F}{r_jp_j} \leq \E[T(j,F)]$$
for any $j \in \cN$, where the second inequality is from \eqref{eqn:lower_upper_bound}. Hence $i^*$ is the static optimal channel, that is, $i^* = i_{so}(F)$. This establishes (b), completing the proof.
\end{IEEEproof} 

\smallskip
Outside of Proposition \ref{prop:threshold_optimal_policy}, however, there are many instances where the max-throughput channel is not static optimal and other channels can perform better for smaller file sizes. In such cases, we would like to discuss how much time the static optimal policy can save against the max-throughput policy.

\begin{corollary}\label{corollary:3.3}
Let $m_i \triangleq p_i/p_{i^*}$ for $i\in \cN \setminus \{i^*\}$. Consider a file of size $F \in (k\Delta r_{i^*}, (k+1)\Delta r_{i^*})$ for some $k \in \N$. Then, we have
%-------------------
\begin{equation}\label{eqn:lower_bound_static_general}
\begin{split}
   \frac{\E[T(i_{so},F)]}{\E[T(i^*,F)]} 
   \!\leq \!\! \min_{\substack{i\in \cN \setminus \{i^*\}}}\!\! \left\{\!1,\!\frac{F/\Delta r_ip_i + (1\!-\!p_i)/p_i}{(k\!+\!1)(m_i/p_i\!-\!1)}\!\right\}\!.
\end{split}
\end{equation}
%-------------------
\end{corollary}
\begin{IEEEproof}
For the file of size $F \in (k\Delta r_{i^*},(k+1)\Delta r_{i^*})$ with $k \in \N$, we have $r_{i^*}p_{i^*} > r_ip_i$ and $m_i \triangleq p_i/p_{i^*}$ for $i \in \cN \setminus \{i^*\}$. From \eqref{eqn:lower_upper_bound} in the proof of Proposition \ref{prop:threshold_optimal_policy}, we have $\E[T(i,F)] \leq F/r_ip_i+\Delta(1-p_i)/p_i$.
Moreover, from \eqref{eqn:time_static_policy} we have
\begin{equation}\label{eqn:lower_max-throughput_policy}
    \E[T(i^*,F)] \!>\! \Delta (k+1)(1/p_{i^*}\!-\!1) \!=\!  \Delta (k+1)(m_i/p_{i}\!-\!1).
\end{equation}
Therefore, the upper bound of the time ratio between channel $i$ and channel $i^*$ is shown as follows:
%-------------------
\begin{equation}\label{eqn:lower_bound_static_example}
\begin{split}
   \frac{\E[T(i,F)]}{\E[T(i^*,F)]} &< \frac{F/\Delta r_ip_i + (1-p_i)/p_i}{(k+1)(m_i/p_i-1)}.
\end{split}
\end{equation}
%-------------------
By definition of the static optimal channel and $\E[T(i_{so},F)]\!\leq\! \E[T(i^*,F)]$, we can get the result \eqref{eqn:lower_bound_static_general} by lower bounding \eqref{eqn:lower_bound_static_example} for channel $i \in \cN \setminus \{i^*\}$.
\end{IEEEproof}

\smallskip
Note that by definition $m_i/p_i = 1/p_{i^*} > 1$ so that the upper bound on the ratio $\E[T(i_{so},F)]/\E[T(i^*,F)]$ in \eqref{eqn:lower_bound_static_general} is always in the interval $(0,1]$. Moreover, smaller ratio means better performance of the static optimal policy against the max-throughput policy. To gauge how the parameters of the max-throughput channel could affect the performance of the static optimal policy, suppose we fix $F,\Delta,r_i,p_i$ for all $i\!\in\! \cN\! \setminus\! \{i^*\}$ and the maximum throughput $r_{i^*} p_{i^*}$, while treating $p_{i^*}$ as a variable. The upper bound in \eqref{eqn:lower_bound_static_general} is then monotonically decreasing in $m_i = p_i/p_{i^*}$, and can even approach to $0$ if at least one of $m_i$ is really large, resulting in the huge performance gain of the static optimal channel compared to that of the max-throughput channel. This implies that accessing channel $i^*$ can take much longer time to transmit a file than other channels if its available probability $p_{i^*}$ is very small, which is common in outdoor networks where the max-throughput channel $i^*$ has very high rate but with low available probability  \cite{bicket2005bit}.

Our static optimal policy shows better performance against the max-throughput policy for small files and small $p_{i^*}$. Since the static optimal channel depends on the file size, choosing channels \emph{dynamically} according to its remaining file size can further reduce the expected transfer time. We next formulate the file transfer problem as an instance of the stochastic shortest path (SSP) problem and analyze the performance of the dynamic optimal policy. 

%% file: Sections/MDP_formulation.tex
\section{Dynamic Optimal Policy}\label{mdp formulation and dynamic policy}

Now that we have analyzed the static policies, we turn our attention to feasible dynamic policies for our file transfer problem. We start by first formulating the file transfer problem as a stochastic shortest path (SSP) problem, in which the agent acts dynamically according to the stochastic environment to reach the predefined destination as soon as possible.
Then, we translate this SSP problem into an equivalent shortest path problem, which helps us derive the closed-form expression of the expected transfer time for any given dynamic policy, and we utilize this to obtain the performance analysis of the dynamic optimal policy against the max-throughput policy.

\subsection{Stochastic Shortest Path Formulation}\label{MDP setting}
The SSP problem is a special case of the infinite horizon Markov decision process \cite{bertsekas2019reinforcement}. To make this section self-contained, we explain our problem as a SSP problem.

\textit{State Space} and \textit{Action Space:} We define the state $s \in \cS \triangleq \R_+$ of our SSP as the remaining file size yet to be transmitted. The action $i \in \cN$ is the channel chosen to be sensed at the beginning of each time interval. The objective of our problem is to take the optimal action at each state $s$ which minimizes the expected time to transmit the file of size $F$.

\textit{State Transition:} Denote by $P_{s,s'}(i)$ the transition probability that the SU moves to state $s'$ after taking action $i$ at state $s$. From any given state $s \in \cS \!\setminus\! \{0\}$, the next state under any action $i\in \cN$ depends on the availability of channel $i$. Since the channel is available or unavailable according to an \textit{i.i.d} (over time) Bernoulli distribution, the next state is either the same as the current one if channel $i$ is unavailable, i.e. $P_{s,s}(i) = 1-p_i$; or the next state is $(s\!-\!\Delta r_i)^+ \triangleq \max\{0,s\!-\!\Delta r_i\}$ if channel $i$ is available, i.e. $P_{s,(s\!-\!\Delta r_i)^+}(i) = p_i$. State $0$ is a termination state since there is no file transmission remaining.

\textit{Cost Function:} The cost $c(s,i,s')$ is the amount of time spent in transition from state $s$ to $s'$ after sensing channel $i$. Since the SU can only sense channels at intervals of size $\Delta$, sensing an unavailable channel costs a $\Delta$ second waiting period until the SU can sense next, that is, $c(s,i,s) = \Delta$ for all $s\in \cS \!\setminus\! \{0\},i \in \cN$. Similarly, if the sensed channel is available, the time spent in transmitting is also $\Delta$ seconds, unless the SU finishes transmitting the file early. In the latter case the cost of transmission is $c(s,i,0) = s/r_i$. Overall, the cost of a successful transmission can be written as $c(s,i,(s\!-\!\Delta r_i)^+) = \min\{\Delta,s/r_i\}$ for all $s\in \cS \!\setminus\! \{0\},i \in \cN$. Once the remaining file size reduces to $0$, the SU will end this file transmission immediately with no additional cost incurred, so that $c(0,i,0) \!=\! 0$ for any $i \in \cN$.

Table \ref{tab:ssp_setting} summarizes the state transition and cost function for our file transfer problem. All other cases except the two cases in Table \ref{tab:ssp_setting} have zero transition probability and zero cost.
\begin{table}[!ht]\vspace{-2mm}  
    \centering
    \caption{SSP setting of our file transfer problem}\vspace{-4mm} 
    \label{tab:ssp_setting}
    \begin{tabular}{|c|c|c|c|c|c|}
         \hline
         current state & action& next state & transition & cost \\
         \hline
         $s>0$ & $i$ & $s$ & $1-p_i$ & $\Delta$ \\
         \hline
         $s>0$ & $i$ & $(s\!-\!\Delta r_i)^+$ & $p_i$ & $\min\{\Delta,s/r_i\}$ \\
         \hline
    \end{tabular}
    \vspace{-2mm}
\end{table}

Our dynamic policy\footnote{There always exists an optimal policy $\pi^*$ to be deterministic in the SSP problem, as proved in Proposition 4.2.4 \cite{bertsekas2019reinforcement}. Thus, we restrict ourselves to the class of deterministic policies in this paper.} is written as a mapping $\pi : S \to \cN$, where $\pi(s)\in\cN$ denotes the channel chosen for sensing when the current state (remaining file size) is $s$. For any policy $\pi$, we have $T(\pi,0) = 0$ at the termination state. Our goal in this SSP problem is to find the dynamic optimal policy $\pi^*(F)$ that minimizes the expected transfer time for the file size $F$, which can be derived from a variety of methods such as value iteration, policy iteration and dynamic programming \cite{bertsekas2019reinforcement}.

\subsection{Performance Analysis}\label{performance analysis}
\begin{figure}[H]
    \vspace{-2mm}
    \centering
    \subfloat{\includegraphics[width=\columnwidth]{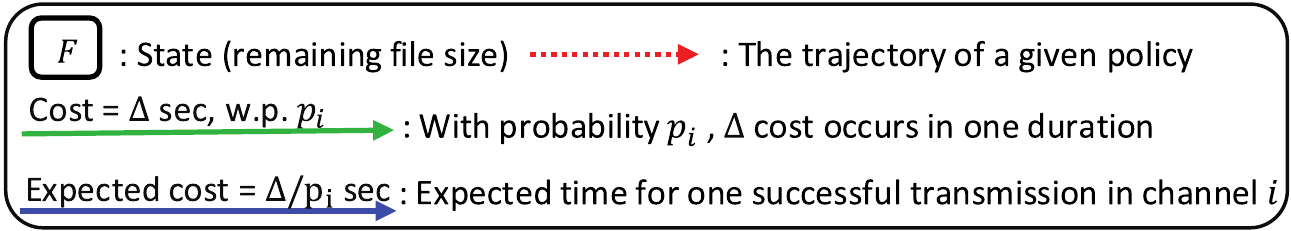}}\\
    \subfloat[(a) Transformation of the SSP problem.\label{fig:deterministic_figure}]{\includegraphics[,width=0.3\columnwidth]{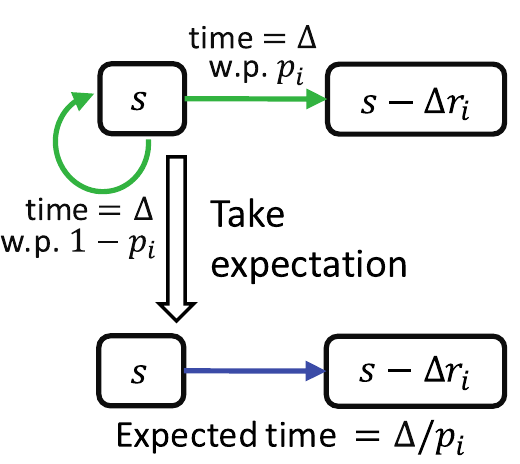}}\quad
    \subfloat[(b) Equivalent shortest path diagram.\label{fig:mdp_instance}]{\includegraphics[width=0.63\columnwidth]{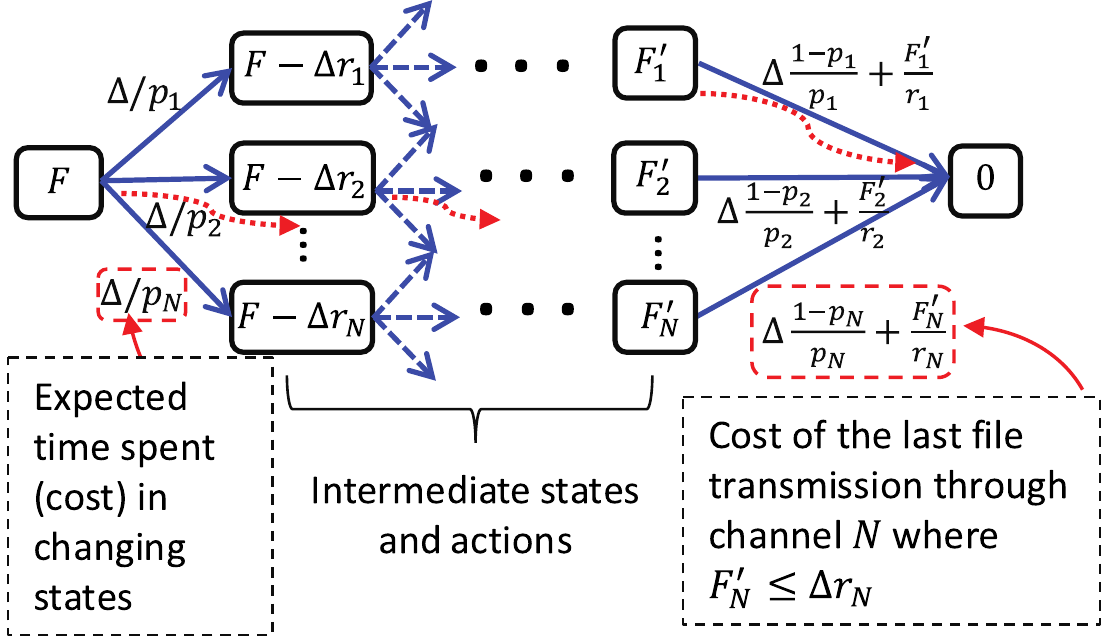}}
    \caption{Illustration to translate the file transfer problem into an equivalent shortest path problem.}
     \label{fig:equivalent_shortest_path}\vspace{-2mm}
\end{figure}

For ease of exposition, we introduce additional notations here. By a \emph{successful transmission}, we refer to state transitions of the form $s \to s'$. This is denoted by the horizontal green line in Figure \ref{fig:equivalent_shortest_path}(a) connecting states $s$ and $s' \triangleq s - \Delta r_i > 0$, and should be distinguished from the self-loop $s \to s$, which implies the sensed channel was unavailable. As shown in Figure \ref{fig:equivalent_shortest_path}(a), taking expectation helps get rid of these self-loops by casting the original SSP to a deterministic shortest path problem in expectation. The cost associated with each link is then the expected time it takes to transit between the states. Figure \ref{fig:equivalent_shortest_path}(b) shows the underlying network for the shortest path problem, where each link is a channel chosen to be sensed and each path from \emph{source} $F$ to \emph{destination} $0$ corresponds to a policy $\pi \in \Pi(F)$. The path-length or the number of links traversed from $F$ to $0$ under any given policy $\pi$ then becomes the total number of successful transmissions needed by that policy to complete the file transfer, which we denote by $|\pi|$.
%----------------------------
%----------------------------

For any policy $\pi \in \Pi(F)$ and $n\in\{1,\cdots,|\pi|\}$, let $F_n$ denote the remaining file size right before the $n$-th successful transmission. Then for all $n\in \{2, \cdots, |\pi|\}$, we have the recursive relationship: $F_n = F_{n-1} - \Delta r_{\pi(F_{n-1})}$, starting with $F_1 \!=\! F$ and ending with $F_{|\pi|+1}=0$. Given a file size $F$, each policy $\pi \in \Pi(F)$ can then be written in a vector form as $\pibf = [\pi(F_1), \pi(F_2), \dots \pi(F_{|\pi|})]^T$.
With this in mind, we can derive a closed-form expression of the expected transfer time for any dynamic policy in the following proposition.

%----------------------------
%--------PROPOSITION---------
%----------------------------
\begin{proposition}\label{prop:expected_time_dynamic_policy}
    Given a file of size $F$, the expected transfer time of a dynamic policy $\pi$ is written as
%----------------------------
\begin{equation}\label{eqn:recursion_closed_form}
\begin{split}
       \E[T(\pi,F)]  \!=\! \Delta \!\!\left(\!\sum_{n=1}^{|\pi|\!-\!1}\! \frac{1}{p_{\pi(F_{n})}} \!\!+\!\!\frac{1\!-\!p_{\pi(F_{|\pi|})}}{p_{\pi(F_{|\pi|})}}\!\!\right)\!\!+\!\! \frac{F_{|\pi|}}{r_{\pi(F_{|\pi|})}}. 
\end{split}
\end{equation}
\end{proposition}
\begin{IEEEproof}
With our notation $E[T(\pi,s)]$ in mind, the Bellman equation for any fixed policy $\pi \in \Pi(F)$ (Proposition 4.2.3 in \cite{bertsekas2019reinforcement}) is shown as
%----------------------------
\begin{equation}\label{eqn:Bellman_policy}
    E[T(\pi,s)] \!=\! \sum_{s' \in \cS}\! P_{s,s'}(\pi(s))\left[c\left(s,\pi(s),s'\right) + E[T(\pi,s')]\right].
\end{equation}
%----------------------------
The transition probability and cost function in section \ref{MDP setting} are defined as $P_{s,s}(i) = 1-p_i, P_{s,(s-\Delta r_i)^+}(i) = p_i, \forall s \in \mathcal{S} \!\setminus\! \{0\}, i \in \cN$, and
\begin{equation*}
    c\left(s,i,(s-\Delta r_i)^+\right) = \begin{cases} \min \{\Delta, s/r_i\} & s \in \mathcal{S} \!\setminus\! \{0\}, i \in \cN \\ 0 & s = 0, i \in \cN.  \end{cases}
\end{equation*}
Then, by substituting $P_{s,s'}(i)$ and $c(s,i,s')$ with our transition probability and cost function defined above, \eqref{eqn:Bellman_policy} can be written as
%----------------------------
\begin{equation*}
\begin{split}
     \E[T(\pi,F)] =& (1-p_{\pi(s)}) \left(\Delta + \E[T(\pi,F)]\right) \\ 
     &+ p_{\pi(s)} \left(\min\left\{\Delta ,\frac{s}{r_{\pi(s)}}\right\} + \E[T(\pi,F_2)]\right).
\end{split}
\end{equation*}
%----------------------------

Recall that $F_n$ is the remaining file size right before the $n$-th successful file transmission given a policy $\pi$ and $F_{n+1} =F_n - \Delta r_{\pi(F_n)}$ for $n = 1,2,\cdots,|\pi|-1$, we can generalize this recursion to two adjacent states $F_n,F_{n+1} \in \cS$ in policy $\pi$ such that
%----------------------------
\begin{equation}\label{eqn:recursion}
\begin{split}
  \E[T(\pi,F_n)] = &\Delta\frac{1-p_{\pi(F_n)}}{p_{\pi(F_n)}} + \min\left\{\Delta ,\frac{F_n}{r_{\pi(F_n)}}\right\} \\ &+   \E[T(\pi,F_{n+1})].    
\end{split}
\end{equation}
%----------------------------
When $n = 1,2,\cdots,|\pi|-1$, the player fully spends $\Delta$ time in each successful transmission and the file transfer task has not been done yet $(F_n> \Delta r_{n})$, so that $\min\left\{\Delta ,F_n/r_{n}\right\} = \Delta$. For the last successful transmission $n=|\pi|$, we have 
$$\E[T(\pi,F_{|\pi|})] = \Delta\frac{1-p_{\pi_{|\pi|}}}{p_{\pi_{|\pi|}}} + \frac{F_{|\pi|}}{r_{\pi_{|\pi|}}}.$$ 
Thereby recursively solving the above equation gives \eqref{eqn:recursion_closed_form}.
\end{IEEEproof}
%----------------------------

\smallskip
In \eqref{eqn:recursion_closed_form}, the first summation is the cumulative expected transmission time, or the cost, to the $(|\pi|-1)$-th successful transmission, with the last two terms being the expected transmission time of the last successful transmission. Proposition \ref{prop:expected_time_dynamic_policy} also includes the expected transfer time of the static policy as a special case. Recall that $k_{i} = \floor{F/\Delta r_i}$ and $\alpha_{i} = F/\Delta r_i - k_i$ in Proposition \ref{prop:expected_time_static}. When applied to a static policy for any channel $i \in \cN$, we have $|\pi| = k_i + \Char_{\{\alpha_i > 0\}}$ and $p_{\pi(F_n)} = p_i$ for all $n = 1,2,\cdots,|\pi|$. The recursive relationship becomes: $F_n = F_{n-1} - \Delta r_i$, implying that $F_n = F - (n-1)\Delta r_i$. Then, we have $F_{|\pi|} = F - (|\pi|-1)\Delta r_i = \alpha_i\Delta r_i$ if $\alpha_i > 0$. Otherwise, $F_{|\pi|} = \Delta r_i$. Substituting these into \eqref{eqn:recursion_closed_form} gets us \eqref{eqn:time_static_policy}.

The common folklore around the max-throughput policy is that it would lead to the minimal file transfer time of $F/r_{i^*}p_{i^*}$. Our next result shows this is too optimistic and not achieved in general even under the dynamic optimal policy. 

\begin{proposition}\label{prop:lower_bound_coincide}
For any file size $F$ and any dynamic policy $\pi \in \Pi(F)$, we have $\E[T(\pi,F)] \geq \E[T(\pi^*,F)] \geq F/r_{i^*}p_{i^*}$.
Moreover, $i^* = \pi^*(F)$ for $F = k\Delta r_{i^*}, k \in \Z_+$.
\end{proposition}
\begin{IEEEproof}
From \eqref{eqn:recursion_closed_form} we have
\begin{equation*}
\begin{split}
   &\E[T(\pi,F)] \\ =& \sum_{n=1}^{|\pi|\!-\!1} \frac{\Delta r_{\pi(F_n)}}{r_{\pi(F_n)}p_{\pi(F_n)}} \!+\! \frac{\Delta\!\left(1\!-\!p_{\pi\left(F_{|\pi|}\right)}\right)\!r_{\pi\left(F_{|\pi|}\right)} \!+\! F_{|\pi|}p_{\pi\left(F_{|\pi|}\right)}}{r_{\pi\left(F_{|\pi|}\right)}p_{\pi\left(F_{|\pi|}\right)}} \\
   \geq& \sum_{n=1}^{|\pi|-1} \frac{\Delta r_{\pi(F_n)}}{r_{i^*}p_{i^*}} +\frac{\Delta\left(1\!-\!p_{\pi\left(F_{|\pi|}\right)}\right)r_{\pi\left(F_{|\pi|}\right)} \!+\! F_{|\pi|}p_{\pi\left(F_{|\pi|}\right)}}{r_{i^*}p_{i^*}} \\
   \geq& \frac{1}{r_{i^*}p_{i^*}}\left(\sum_{n=1}^{|\pi|-1}\Delta r_{\pi(F_n)} + F_{|\pi|}\right) = \frac{F}{r_{i^*}p_{i^*}},
\end{split}
\end{equation*}
where the first inequality comes from the fact that $r_{i^*}p_{i^*} \geq r_i p_i$ for all $i \in \cN$. The second inequality is from our definition of $F_{|\pi|}$ which implies that $F_{|\pi|} \leq \Delta r_{\pi_{|\pi|}}$.

When file size $F$ is an integer multiple of $\Delta r_{i^*}$, i.e., $F=k \Delta r_{i^*}$ for some $k \in \Z_+$, we have $\E[T(i^*,F)] = F/r_{i^*}p_{i^*}\leq \E[T(\pi^*,F)] \leq E[T(\pi,F)]$. Since $\E[T(i^*,F)] \geq \E[T(\pi^*,F)]$, we have $\E[T(\pi^*,F)] = F/r_{i^*}p_{i^*}$, and the max-throughput policy coincides with the dynamic optimal policy.
\end{IEEEproof}

\smallskip
As shown in \eqref{eqn:static_lower_bound}, $F/r_{i^*}p_{i^*}$ is always the lower bound on the transfer time for any static policy. Proposition \ref{prop:lower_bound_coincide} strengthens this by showing that the same is true even for the dynamic optimal policy. Similar to condition (b) in Proposition \ref{prop:threshold_optimal_policy} for the static optimal policy, the dynamic optimal policy $\pi^*(F)$ also coincides with the max-throughput policy $i^*$ when the file size is an integer multiple of $\Delta r_{i^*}$, while we no longer have the finite threshold $H$ as in Proposition \ref{prop:threshold_optimal_policy}(a). We next give bounds to quantify the performance of the dynamic optimal policy with respect to the max-throughput policy. 
%----------------------------
%--------PROPOSITION---------
%----------------------------
\begin{corollary}\label{corollary:dynamic_lower_upper_bound}
Let $m_i \triangleq p_i/p_{i^*}$ for $i\in \cN \setminus \{i^*\}$. Consider a file of size $F \in (k\Delta r_{i^*}, (k+1)\Delta r_{i^*})$ for some $k \in \N$. Then, we have
\begin{equation*}
\begin{split}
    &\frac{F}{F+\Delta \Char_{\{\alpha_i > 0\}}(1-p_{i^*})r_{i^*}}\leq \frac{\E[T(\pi^*,F)]}{\E[T(i^*,F)]}  \\
    &\leq \min_{\substack{i \in \cN \setminus \{i^*\}}}\left\{1,\frac{\frac{F}{\Delta r_i} + 1-p_i - k m_i\frac{r_{i^*}p_{i^*}-r_ip_i}{r_ip_i}}{(k+1)(m_i-p_i)}\right\}.
\end{split}
\end{equation*}
\end{corollary}
\begin{IEEEproof}
We first define a suboptimal policy $\pi_{heu}$ and then use $E[T(\pi^*,F)] \leq E[T(\pi_{heu},F)]$ for our proof.
The suboptimal policy $\pi_{heu}$ is defined as sensing and accessing the max-throughput channel $i^*$ to transmit the file of size $F_1 = n\Delta r_{i^*}$, then following a static optimal policy for the remaining file of size $F_2 = F - n\Delta r_{i^*}$. $n$ is an integer chosen from $0$ to $k=\floor{F/\Delta r_{i^*}}$.
Then, the expected transfer time of the dynamic optimal policy is always smaller than that of the dynamic suboptimal policy, that is, 
\begin{equation*}
    \begin{split}
        \E[T(\pi^*,F)] &\leq \E[T(\pi_{heu},F)]\\
        &=\E[T(i^*,n\Delta r_{i^*})]+ \E[T(i_{so},F-n\Delta r_{i^*})]\\
        &\leq \Delta n/p_{i^*} +\min_{\substack{i \in \cN \!\setminus\!\{i^*\}}} \left\{ \frac{F-n\Delta r_{i^*}}{r_ip_i} +  \Delta\frac{1-p_i}{p_i}\right\} \\
        &= \!\!\min_{\substack{i \in \cN \!\setminus\!\{i^*\} }}\!\!\left\{\!\frac{F}{r_ip_i} + \Delta\left( \frac{1-p_i}{p_i} \!-\! n\frac{r_{i^*}p_{i^*}\!-\!r_ip_i}{r_ip_ip_{i^*}}\right)\!\right\}\!,
    \end{split}
\end{equation*}
where the second inequality comes from \eqref{eqn:time_static_policy}, and \eqref{eqn:lower_upper_bound} in the proof of Proposition \ref{prop:threshold_optimal_policy}. It shows monotonically decreasing in $n$ such that we can choose $n=k$ to get the smallest upper bound for $\E[T(\pi^*,F)]$. Moreover, we have $E[T(i^*,F)] > \Delta (k+1)(m_i/p_{i}-1)$ from \eqref{eqn:lower_max-throughput_policy}. Hence, the upper bound of the ratio $\E[T(\pi^*,F)]/\E[T(i^*,F)]$ in Corollary \ref{corollary:dynamic_lower_upper_bound} is proved.

For the lower bound of the ratio, by using Proposition \ref{prop:lower_bound_coincide} we have $E[T(\pi^*,F)] \geq F/r_{i^*}p_{i^*}$. Together with \eqref{eqn:static_lower_bound}, we have
\begin{equation*}
\begin{split}
     \frac{\E[T(\pi^*,F)]}{\E[T(i^*,F)]} &\geq \frac{F/r_{i^*}p_{i^*}}{F/r_{i^*}p_{i^*} + \Delta \Char_{\{\alpha_{i^*}>0\}}  (1\!-\!\alpha_{i^*}) (1\!-\!p_{i^*})/p_{i^*}} \\
     &\geq \frac{1}{1 + \Delta \Char_{\{\alpha_{i^*}>0\}}(1-p_{i^*})r_{i^*}/F},
\end{split}
\end{equation*}
where the second inequality comes from $1-\alpha_{i^*}\leq 1$. This completes the proof.
\end{IEEEproof}

\smallskip
To better understand Corollary \ref{corollary:dynamic_lower_upper_bound}, we analyze how the parameters of the max-throughput channel could impact the performance of the dynamic optimal policy. Similar to Corollary \ref{corollary:3.3}, small value of ${\E[T(\pi^*,F)]}/{\E[T(i^*,F)]}$ implies that the dynamic optimal policy offers significant saving in time over the max-throughput policy. We note that Corollary \ref{corollary:dynamic_lower_upper_bound} tightens the upper bound with an extra negative term in the numerator, compared to that in Corollary \ref{corollary:3.3}, potentially providing greater savings in time as we extend the policy from static optimal to dynamic optimal.

In contrast to Proposition \ref{prop:threshold_optimal_policy} that max-throughput policy is good enough for $F \geq H$, Corollary \ref{corollary:dynamic_lower_upper_bound} tells us that there is always some reduction in file transfer time even for large file size $F$ under the dynamic optimal policy. This is because the extra negative term in the numerator can be large, since $k=\floor{F/\Delta r_{i^*}}$ could be big for large $F$, implying that the second argument in the $\min\{\cdot,\cdot\}$ function may no longer be increasing in $F$. Note however that the reduction in transfer time would be minimal for large file sizes since the lower bound in Corollary \ref{corollary:dynamic_lower_upper_bound} will rise to $1$ as $F$ goes to infinity.

%% file: Sections/Practical_concern.tex
\section{Practical considerations}\label{online learning}
While the dynamic optimal policy gives a smaller expected transfer time, the computational cost of solving the shortest path problem is still a concern. Besides, we face a scaling problem when the file size differs from each, effectively changing the underlying ``graph'' in the corresponding shortest path problem. This warrants re-computation of the dynamic optimal polity for each file size, which would be unacceptable in reality. In this section, we discuss a mixed-integer programming formulation and propose a heuristic policy to balance the performance and the computational cost. We also consider the case where the SU doesn't know about the channel parameters beforehand and it must sense and access channels \emph{on the fly} to find the optimal policy.

\subsection{Mixed-Integer Programming Formulation}\label{subsection:mix-integer}
Dynamic programming problems often have equivalent integer or mixed integer programming formulations as well \cite{della2017exact}. For our shortest path problem, however, we can leverage the fact that the cost of an action is the same for each state (before the last successful transmission) to notably reduce the size of the solution space of the mixed integer formulation, especially for large file sizes where the curse of dimensionality is most felt.

Before putting forward our equivalent mixed integer programming problem, We first provide an alternate expression for \eqref{eqn:recursion_closed_form} as below.

\begin{proposition} \label{prop:equivalent_closed_form}
Given a policy $\pi\in\Pi(F)$ for any file of size $F$, let $x_i \!\triangleq\! \sum_{n=1}^{|\pi|\!-\!1} \Char_{\{\pi(F_n)=i\}}$, and let $y_i \triangleq \Char_{\{ \pi(F_{|\pi|})=i\}} F_{|\pi|}/\Delta r_{\pi(F_{|\pi|})}$. Then, we can rewrite \eqref{eqn:recursion_closed_form} as 
%================
\begin{equation*}
\begin{split}
    E[T(\pi,F)] \!=& \Delta \!\sum_{i=1}^N \Bigg[ \sum_{n=1}^{|\pi|-1} \frac{\Char_{\{\pi(F_n)=i\}}}{p_i} \\
    &+ \Char_{\{\pi(F_{|\pi|})=i\}} \!\left( \!\frac{1\!-\!p_{\pi(F_{|\pi|})}}{p_{\pi(F_{|\pi|})}} \!+\! \frac{F_{|\pi|}}{\Delta r_{\pi(F_{|\pi|})}} \!\right) \!\Bigg] \\
    =& \Delta \sum_{i=1}^N\frac{x_i}{p_i} + \Char_{\{y_i>0\}} \frac{1-p_i}{p_i} + y_i. ~~~~~~~~~~~~~~\qedsymbol
\end{split}
\end{equation*}
%=================
\end{proposition}
In the above, $x_i$ counts the total number of successful transmissions through any channel $i\in\cN$ except for the last ($|\pi|$-th) transmission. On the other hand, $y_i$ is zero if channel $i$ is not the last one sensed under policy $\pi$, or else is equal to the fraction of the $\Delta$ second interval used for the $|\pi|$-th successful transmission.

A combination of $\vx=[x_i]_{i\in\cN}$ and $\vy=[y_i]_{i\in\cN}$ can be the same for multiple optimal policies, which have the same expected time although different order in which the channels are sensed. This means there can be $(\ones^T\vx)!/\prod_{i \in \cN} (x_i!)$ policies with the same expected transfer time - a number which can be really large for large $F$. Condensing the state space by preventing the solver from considering these $(\ones^T\vx)!/\prod_{i \in \cN} (x_i!)$ many policies individually can significantly reduce computation time. This can be done by solving the following mixed-integer programming problem over the set of all $\vx,\vy$ that correspond to feasible dynamic policies.
%=======
\begin{equation}\label{mixed-integer}
    \begin{split}
        \min_{\vx,\vy} \ \ \  &\Delta  \!\sum_{i=1}^N \left[\frac{x_i}{p_i} + \Char_{\{y_i>0\}}\frac{1-p_i}{p_i} + y_i\right] \\
        s.t. \ \ \  &\Delta \!\sum_{i=1}^{N}\left(x_i r_i + y_i r_i\right) = F,\\
        &y_iy_j = 0, \ \ \ \ \ \ \ \ \ \ \ \ \ \ \ \ \ \ \ \forall i\neq j \in \cN,\\
        &x_i \in \Z_+,\ y_i \in [0,1), \ \ \forall i \in \cN.
    \end{split}
\end{equation}
%=======

In the above optimization problem, the first constraint ensures that the choice of $\vx,\vy$ guarantees the transmission of the entire file by adding up to $F$ when multiplied by $\Delta$ and the respective channel rates. The second constraint makes sure that at most one of the $y_i$'s is positive to ensure that the last transmission, if any, is assigned only to one channel. Once an optimal solution $\vx^* = [x_i^*]$ and $\vy^*=[y_i^*]$ is obtained, we can construct a corresponding dynamic optimal policy $\pi^*$ by setting the sole channel $i$ for which $y_i>0$ as the last channel for transmission under the policy, and the first $|\pi^*|-1$ transmissions can be according to any permutation of assignments\footnote{For example, if $\vx=[0,3,2]^T$ for a three channel system, then the $5$ successful transmissions will be $3$ of them via channel 2 and $2$ of them via channel $3$, in no particular order.} from the vector $\vx$. In practice, we usually consider transmitting the file in channel $i$ successfully for $x_i$ counts and then switching to the closest channel $j$ for $x_j$ successful transmissions, which can help reduce the switching delay and energy cost. In this way, we expect that searching for a dynamic optimal policy over this condensed space would be much quicker than searching over the entire set of paths for the shortest path. 

Now, we formally compare the computational complexity of Dijkstra algorithm and mixed integer programming, which is known to be NP-hard in general \cite{conforti2014integer}. In the following lemma, denote by $r_{\min}$ the minimum channel rate over all channels and $\ceil{\cdot}$ the ceiling function, we show that Dijkstra algorithm for the underlying graph in Figure \ref{fig:equivalent_shortest_path} is also an NP-hard problem w.r.t the file size $F$.
\begin{lemma}\label{lemma:complexity}
The worst-case computational complexity of Dijkstra algorithm is $O(\ceil{F/\Delta r_{\min}} N^{\ceil{F/\Delta r_{\min}}}\log(N))$.
\end{lemma}
\begin{IEEEproof}
Dijkstra algorithm is known to have the time complexity $O(E\log(V))$ \cite{bertsekas2012dynamic}, where $E$ is the number of edges and $V$ is the number of nodes of the underlying graph. We now consider the tree-like underlying graph structure, $F$ being the source node and $0$ being the destination, as the worst case. The tree-depth is $\ceil{F/\Delta r_{\min}}$ and the node in each level contains $N$ child nodes, resulting in $O(N^{\ceil{F/\Delta r_{\min}}})$ edges in total. In this case, $E = O(N^{\ceil{F/\Delta r_{\min}}})$, $V \approx E $, and we have
\begin{equation*}
    O(E\log(V)) = O(\ceil{F/\Delta r_{\min}} N^{\ceil{F/\Delta r_{\min}}}\log(N)).
\end{equation*}
\end{IEEEproof}
\smallskip
Lemma \ref{lemma:complexity} shows that the time complexity of Dijkstra algorithm is exponential in the file size $F$. As a result, both Dijkstra algorithm and mixed integer programming are NP-hard w.r.t the file size $F$. However, from the practical implementation, we address that the underlying graph of the shortest path problem is not given upfront. Later in Section \ref{computation_policy_complexity}, we will point out that the SU needs to first generate the graph of the shortest path problem in the real-world implementation, which includes all possible paths from the source node $F$ to the destination $0$ and involves a huge computational overhead, before feeding it into the dynamic programming solver. We observe that generating such graph is the most time-consuming step, while the mixed integer programming solver doesn't need such overhead and runs faster in practice.

\subsection{Performance-Complexity Trade Off}\label{subsection: reusable}
Transmitting different-sized files is very common in the real world. For example, a short text-only email only takes up $5$KB, one five-page paper is around $100$KB and the average size of web page is $2$MB, all implying that the file size may vary greatly \cite{mantuano_2016}. However, due to the nature of the shortest path problem, a change in file size induces a change in the underlying graph. If the goal is to always determine the best solution, the only option is to recompute the dynamic optimal policy for every different file size. This would not be scalable in applications where minimal computation is required, and policies that can be promptly modified and reused across different file sizes with performance guarantees would be highly desirable.

To avoid heavy computation for each file size to obtain the dynamic optimal policy, we propose a heuristic policy that utilizes the max-throughput policy and the static optimal policy in order to reduce the computational cost, while still maintaining considerable performance gain. Note that the max-throughput policy coincides with the dynamic optimal policy when the file size is an integer multiple of $\Delta r_{i^*}$ according to Proposition \ref{prop:lower_bound_coincide}. We also know that the static optimal policy significantly outperforms the max-throughput policy especially for smaller file sizes. Combining these two policies, by transmitting file through max-throughput channel until the remaining file size becomes `small' so as to apply the static optimal policy for the rest,\footnote{We can choose remaining file size to be smaller than $\Delta r_{i^*}$ to apply the static optimal policy.} will strike the right balance between computational complexity and achievable performance gain. From the computational-cost-saving perspective, max-throughput policy is fixed and known to the SU. The closed-form expression $E[T(i,F)]$ for static policy (channel $i$) is also known to the SU. Since our heuristic policy includes the ``min'' function on $\E[T(i,F)]$ for all channel $i\in\mathcal{N}$, its computational complexity is $O(N)$, which is much smaller than that of the dynamic optimal policy whose complexity is polynomial in $N$ as described in Lemma \ref{lemma:complexity}.

With this motivation in mind, we divide file size $F := F_1 + F_2$ into two parts: $F_1 = n \Delta r_{i^*}$ ($n = 0,1,\cdots,\floor{F/\Delta r_{i^*}}$) and $F_2 = F- F_1$. The heuristic policy $\pi_{heu}$ is defined as follows: The SU first transmits the file of size $F_1$ through the max-throughput channel $i^*$ and then sticks to the static optimal policy $i_{so}(F_2)$ for the remaining file of size $F_2$.\footnote{For $F$ being integer multiple of $\Delta r_{i^*}$, we have $F_1 = F$ and $F_2 = 0$, then the heuristic policy coincides with the max-throughput policy, which is also the dynamic optimal policy in view of Proposition \ref{prop:lower_bound_coincide}.} We have shown in the proof of Corollary \ref{corollary:dynamic_lower_upper_bound} that the upper bound of ratio $E[T(\pi_{heu},F)]/E[T(i^*,F)]$ is monotonically decreasing in $n$. Therefore, $n = k$ potentially gives us the smallest upper bound of the ratio (the same upper bound in Corollary \ref{corollary:dynamic_lower_upper_bound}). Moreover, we have explained after Corollary \ref{corollary:dynamic_lower_upper_bound} that the upper bound is smaller than that of the static optimal policy in Corollary \ref{corollary:3.3}. These arguments suggest that the heuristic policy $\pi_{heu}$ with $n=k$ can potentially offer smaller delay than other candidates with different values of $n$. Besides, our heuristic policy can further reduce the computational cost for a set of files sharing the same remaining file size $F_2$, for which the static optimal policy $i_{so}(F_2)$ has already been found and no further re-computation is needed.

\subsection{Unknown Channel Environment}\label{subsection:online_implementation}
We now consider the setting where the SU does not know the available probability $p_i$ for any channel $i\in\cN$ and only knows the rate $r_i$ --- a commonly analysed setting in the OSA literature \cite{mohamedou2017bayesian,dai2012efficient}. The SU has no alternative but to observe the states of these channels when it tries to access them, and build its own estimations of channel probabilities. In this extended setting, we study our problem as an online shortest path problem, which has been widely studied in \cite{gai2012combinatorial,chen2013combinatorial,talebi2017stochastic} for different kinds of cost functions. \cite{talebi2017stochastic} proposed a Kullback-Leibler source routing (\textit{KL-SR}) algorithm to an online routing problem with geometrically distributed delay in each link, which coincides with our link cost in the underlying graph of the shortest path problem shown in Figure \ref{fig:equivalent_shortest_path}(b).

For our purpose, we modify \textit{KL-SR} algorithm; key differences being that we let the file size vary across the episodes, allowing a different underlying graph of the shortest path problem for each episode instead of the fixed underlying graph of the shortest path problem in \cite{talebi2017stochastic}. Algorithm \ref{alg:online_algorithm} describes our online implementation,
where $F^k$ is the file size to be transferred in the $k$-th episode. $n_i(k)$ denotes the number of times channel $i$ has been sensed before the $k$-th episode and $\bar{p}_i(k)$ is the empirical average of channel $i$'s available probability throughout the $k\!-\!1$ episodes so far. With $n_i(k)$ and $\bar{p}_i(k)$, the estimated available probability $\hat{p}_i(k)$ of channel $i$ is then derived from the KL-based index in \cite{talebi2017stochastic}. As mentioned in line $1$ in Algorithm \ref{alg:online_algorithm}, the SU can choose one of the various policies according to which it wishes to perform the file transfer, i.e., dynamic optimal policy $\pi^*$, static optimal policy $i_{so}$, max throughput policy $i^*$ or the heuristic policy $\pi_{heu}$, and then stick to that policy.
Let $\tilde{\E}[T(\pi,F^k)]$ be the estimated expected transfer time of policy $\pi$ at the $k$-th file by using the estimated parameter $\hat{p}_i(k)$ instead of $p_i$ for all $i \in \cN$ in \eqref{eqn:recursion_closed_form}. In line $7$ in Algorithm \ref{alg:online_algorithm}, $\pi^k$ will be computed as $\pi^k \!=\! \argmin_{\pi \in \Pi(F^k)}\tilde{\E}[T(\pi,F^k)]$ for dynamic optimal policy; $\pi^k \!=\! \argmin_{i \in \cN}\tilde{\E}[T(i,F^k)]$ for static optimal policy and $\pi^k \!=\! \argmax r_i \hat{p}_i(k)$ for max-throughput policy. For heuristic policy $\pi_{heu}$,\! $\pi^k\!$ will be computed in the same way as described in Section \ref{subsection: reusable} with estimated parameters $\{\hat{p}_i(k)\}_{i \in \cN}$.
%------------ALGORITHM1----------Online ------
\begin{algorithm}[!h]
	\caption{Online file transfer algorithm}
	\label{alg:online_algorithm}
	\begin{algorithmic}[1]
	    \State Choose the type of policy to use: $\pi^*$, $i_{so}$, $i^*$ or $\pi_{heu}$.
	    \State Apply static policy $i$ in the $i$-th episode to transmit the file of size $F^i$ for $i=1,2,\cdots,N$, and update $n_i(N+1)$, $\bar{p}_i(N+1)$ for $i \in \cN$ at the end of the $N$-th episode 
	    \For {$k \geq N+1$}
	        \State Compute the estimated channel statistics $\{\hat{p}_i(k)\}_{i\in \cN}$ according to the KL-based index in \cite{talebi2017stochastic}.
	        \State Given a file of size $F^k$, compute the policy $\pi^k$ with $\{\hat{p}_i(k)\}_{i \in \cN}$ and observe channel status for the whole file transfer process in this episode. 
	        \State Update $n_i(k+1)$, $\bar{p}_i(k+1)$ for $i \in \cN$.
	    \EndFor
	\end{algorithmic}
\end{algorithm}
%-------------------------------------

The performance of Algorithm \ref{alg:online_algorithm} with varying file sizes (assuming bounded file size) is measured by its \textit{regret} $\E[R(K)]$, which is defined as the cumulative difference of expected transfer time between policy $\pi^k$ at $k$-th file and the targeted optimal policy $\pi_{tar} \in \{\pi^*, i_{so}, i^*, \pi_{heu}\}$ up to the $K$-th file. The regret analysis is nearly the same as Theorem $5.4$ in \cite{talebi2017stochastic}. Let $f(K) = \log(K)+4\log(\log(K))$ and $H = F_{\max}/r_{\min}$, where $F_{\max}$ is the largest possible file size, $r_{\min} \!=\! \min_{i \in [N]} r_i$. Denote by $D = \max_{\pi} \E[T(\pi,F_{\max})]$ the longest expected transfer time and $\epsilon = (1-2^{-\frac{1}{4}})\Delta_{\min}/D$,
\begin{equation*}
    \Delta_{\min} \!=\! \min_{\substack{F \in (0,F_{\max}], \pi\in \Pi(F)/\pi_{tar}}} \E[T(\pi,F)] \!-\! \E[T(\pi_{tar}(F),F)]
\end{equation*}
is the smallest non-zero difference of expected transfer time between any sub-optimal policy $\pi$ and the targeted optimal policy $\pi_{tar}$. Let $p_{\min} \!=\! \min_{i \in [N]} p_i$. The regret bound of Algorithm \ref{alg:online_algorithm} is given in the following theorem.
\begin{theorem}\label{thm:regret_bound}
The gap-dependent regret bound under Algorithm \ref{alg:online_algorithm} is
\begin{equation}\label{eqn:regret_bound}
    \E[R(K)] \!\leq \! \frac{360NH f(K)}{\Delta_{\min}p_{\min}^2} \!+\! 2D\left(4H \!+\!s \sum_{i=1}^n \frac{1}{\epsilon^2 p^2_i}\right).
\end{equation}
\end{theorem}
\begin{IEEEproof}
The proof is nearly the same as the analysis of Theorem 5.4 in Appendix G.B \cite{talebi2017stochastic}. Here we only give the main modifications for our setting.

The first modification comes from the definition of `arm'. In \cite{talebi2017stochastic}, each edge in the graph is treated as a different arm, that is, the status of each edge is observed and estimated separated. In our setting, each edge in the shortest path problem (Figure \ref{fig:equivalent_shortest_path}) represents one of $N$ channels such that each policy (path) may observe one channel multiple times. Then, some summation terms in the proof, previously were over all edges (e.g., (12), (13) in \cite{talebi2017stochastic}), are now over all $N$ channels. 

Second, the \textit{KL-SR} algorithm in \cite{talebi2017stochastic} for dynamic optimal policy works for a fixed source node, which can be interpreted as a fixed file size $F$. Our algorithm deals with the varying file size. Since file size $F_{\max} < \infty$, we only need to change parameter $H$ to be the longest policy length for maximum file size $F_{\max}$ (instead of fixed file size $F$), $\Delta_{\min}$ to be the smallest non-zero difference of expected transfer time between any sub-optimal policy and targeted optimal policy for file size in $(0,F_{\max}]$ (instead of fixed file size $F$) and $D$ to be the longest expected transfer time for file size $F_{\max}$ (instead of fixed file size $F$). Then, the proof will be carried over.
\end{IEEEproof}
The regret \eqref{eqn:regret_bound} scales linearly with the number of channels $N$, instead of the number of edges in the online shortest path problem \cite{talebi2017stochastic}, because each edge in our setting (see Figure \ref{fig:equivalent_shortest_path}) is chosen from one of $N$ channels while each edge in \cite{talebi2017stochastic} is treated as a different `arm'.

%% file: Sections/Simulation.tex
\section{Numerical Results}\label{simulation}
In this section, we present numerical results for file transfer time under four different policies in both offline setting (known $p_i$'s) and online setting (unknown $p_i$'s), using three different channel scenarios as in \cite{bicket2005bit,Srikant_2019}. Through these results, we show the significant time reduction achieved by the dynamic optimal, static optimal and heuristic polices over the max-throughput channel, in line with theoretical analysis.

\subsection{Simulation Setup}
We consider the experimental setup as an IEEE 802.22 system with $8$ different channels. The time duration $\Delta$ is set to $100$ ms, per IEEE 802.22 standard \cite{802.22_standard}. We use three different channel scenarios: \textit{gradual, steep} and \textit{lossy} \cite{bicket2005bit,Srikant_2019}. \textit{Gradual} refers to a case where the available probability of the max-throughput channel is larger than $0.5$. \emph{Steep} is characterized by the available probability of each channel being either very high or very low. \textit{Lossy} means that the available probability of the max-throughput channel is smaller than $0.5$. The channel parameters used for simulation in the above three channel scenarios are given in Table \ref{tab:parameter}. All simulations are run on a PC with AMD Ryzen 1700X and 32G RAM.
\begin{table}[!ht]
    \centering
    \caption{Channel parameters in three channel scenarios}\vspace{-4mm}
     \label{tab:parameter}
    \begin{adjustbox}{width=\columnwidth}
    \begin{tabular}{|c|c|c|c|c|c|c|c|c|}
         \hline
         channel $i$ & $1$ & $2$ & $3$ & $4$ & $5$ & $6$ & $7$ & $8$\\
         \hline
         $r_i$ (Mbps) & $1.5$ & $4.5$ & $6$ & $9$ & $12$ & $18$ & $20$ & $23$ \\
         \hline
         $p_i$ (\emph{gradual}) & $0.95$ & $0.85$ & $0.75$ & $0.65$ & $0.4$ & $0.3$ & $0.2$& $0.1$ \\
         \hline
         $p_i$ (\emph{steep}) & $0.9$ & $0.25$ & $0.2$ & $0.18$ & $0.17$ & $0.16$ & $0.15$ & $0.14$ \\
         \hline
         $p_i$ (\emph{lossy}) & $0.9$& $0.8$ & $0.7$ & $0.4$ & $0.3$ & $0.25$ & $0.2$ & $0.1$ \\
         \hline             
    \end{tabular}
    \end{adjustbox}
    \vspace{-2mm}
\end{table}
%----------------------

\subsection{Computation of Dynamic Optimal Policy}\label{computation_policy_complexity}
To obtain the dynamic optimal policy in both offline and online cases, we utilize the mixed-integer programming formulation as in Section \ref{subsection:mix-integer} to fasten the simulation speed and use the SCIP solver \cite{MaherMiltenbergerPedrosoRehfeldtSchwarzSerrano2016}. For SSP formulation as in Section \ref{MDP setting}, we use policy iteration as a solver. We select $10$ files of sizes in the interval $(0,7]$ (Mb) uniformly at random and compare the total time of computing the dynamic optimal policies of these $10$ files from the mixed-integer programming formulation and SSP formulation. The policy iteration takes $77.94$ seconds, while the SCIP solver only takes $0.23$ seconds. We observe this because any dynamic programming procedure has to effectively first construct the underlying network for the shortest path problem, and then traverse all the possible paths from the source $F$ to destination $0$. This underlying network changes for every different $F$ as well, rendering previous computations useless. However, the mixed-integer programming doesn't need the construction of a network to solve it and SCIP solver only needs to compute one possible solution to \eqref{mixed-integer}. 

\subsection{Online File Transfer Simulation}
%----------------------------
\begin{figure*}[!ht]\vspace{-0mm}
    \centering
    \includegraphics[width=0.97\textwidth]{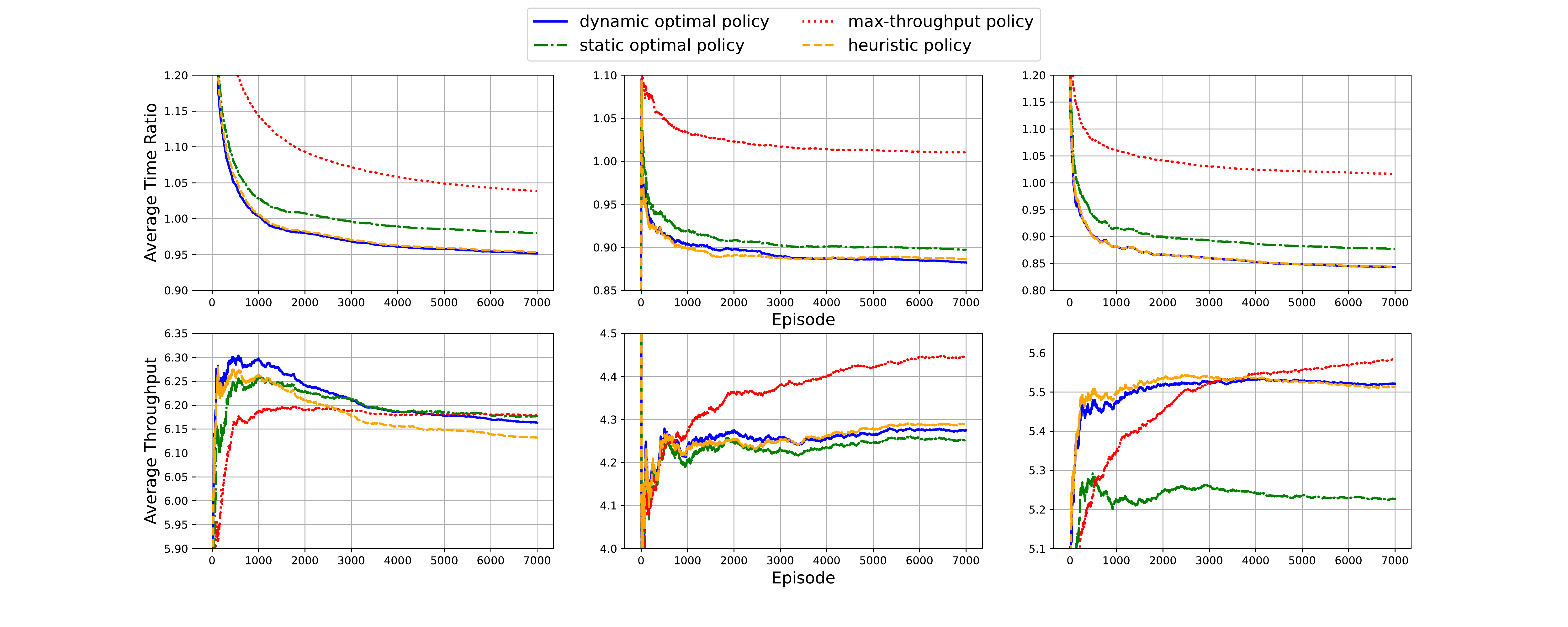}
    \vspace{-4mm}
    \caption{Average time ratio (top row) and average throughput (bottom row) in the online setting. Channel scenarios from left to the right: Gradual; Steep; Lossy.}
    \label{fig:online_simulation_gain}\vspace{-4mm}
\end{figure*}
%----------------------------
\begin{figure*}[!t]
    \vspace{-2mm}
    \centering
    \subfloat[(a) Offline setting, no switching delay. \label{fig:offline_dynamic_vs_max_throughput}]{\includegraphics[width=0.3\textwidth]{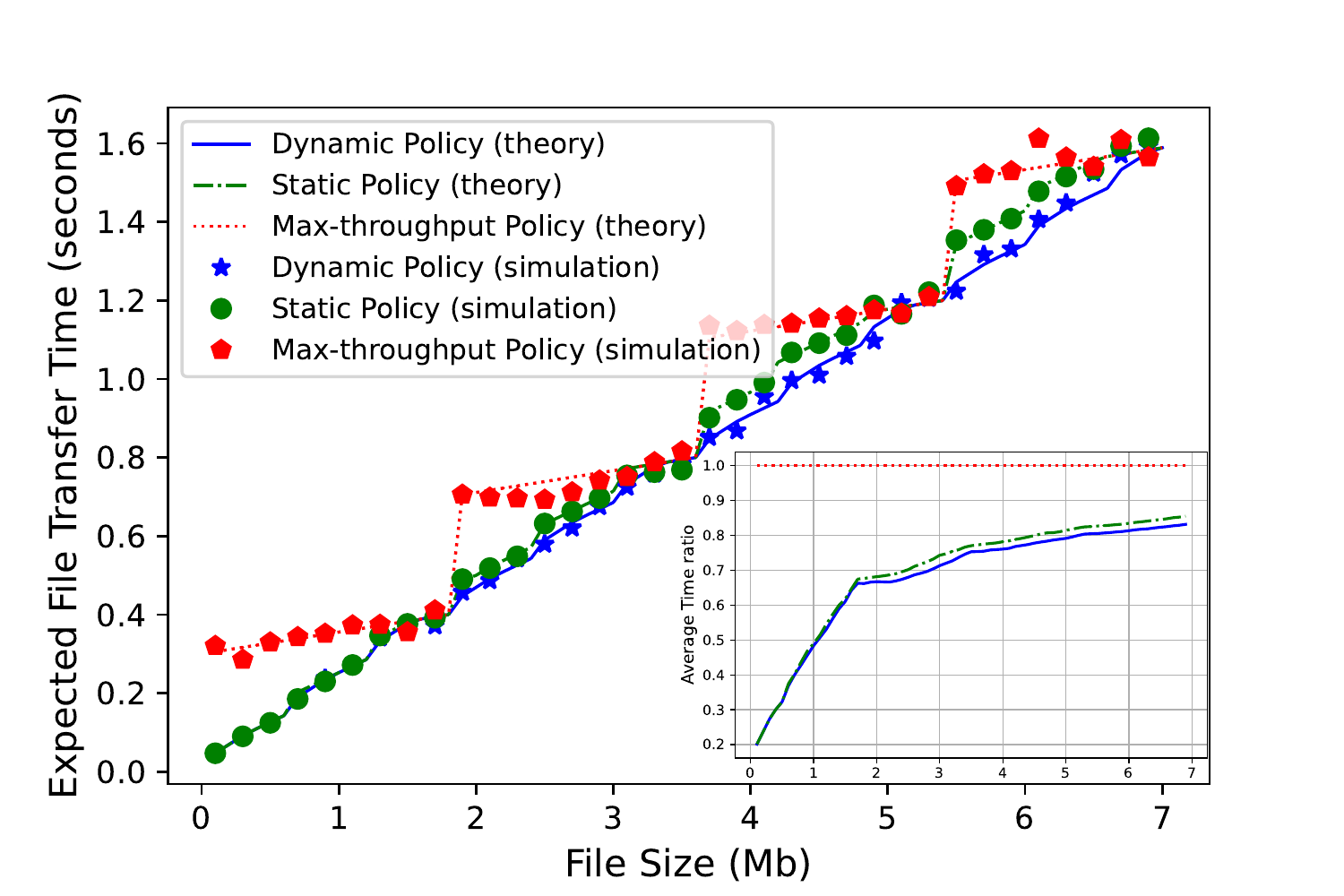}}\quad
    \subfloat[(b) Offline setting, with switching delay. \label{fig:offline_sc}]{\includegraphics[width=0.3\textwidth]{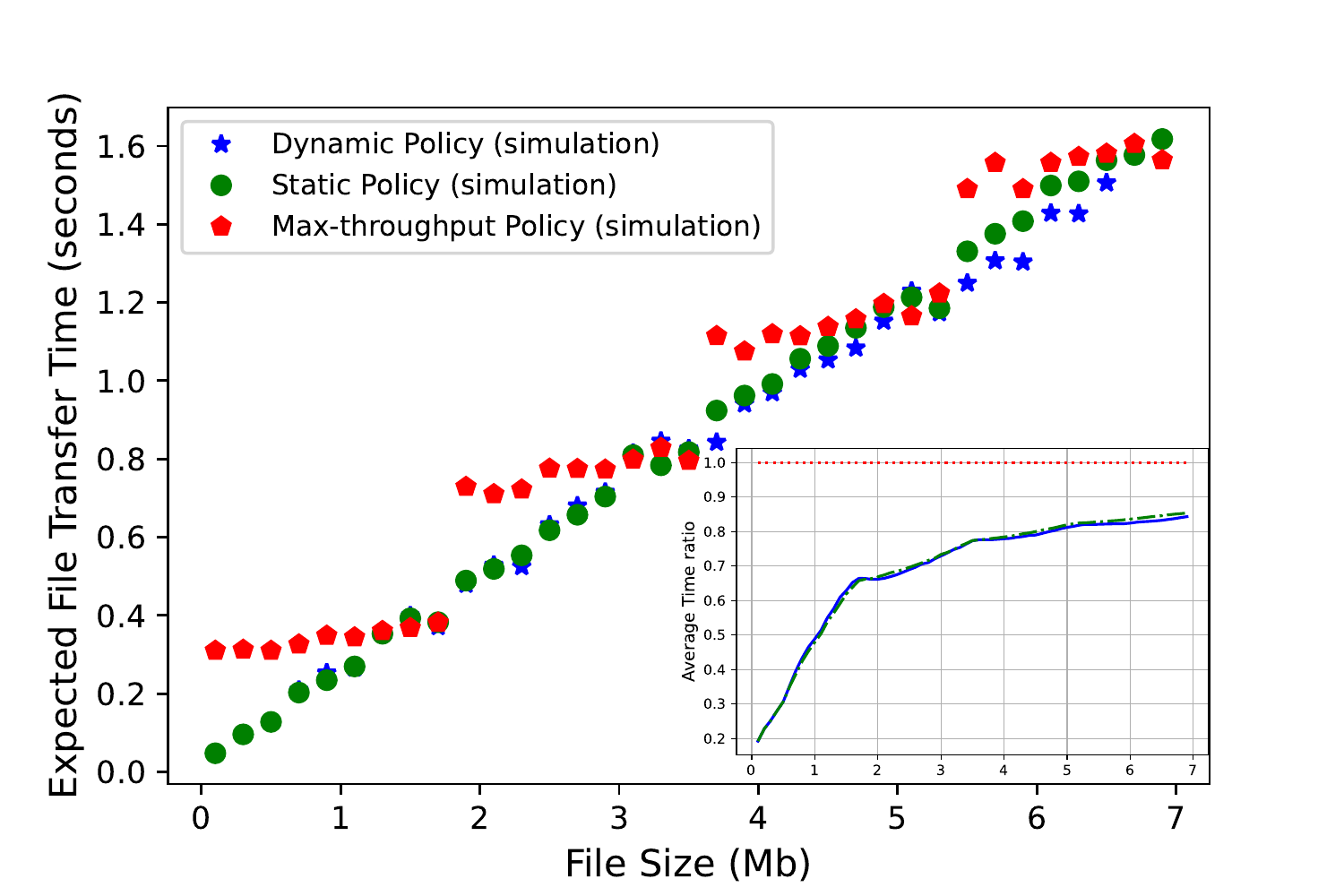}}\quad
     \subfloat[(c) Online setting, with switching delay. \label{fig:online_sc}]{\includegraphics[width=0.335\textwidth]{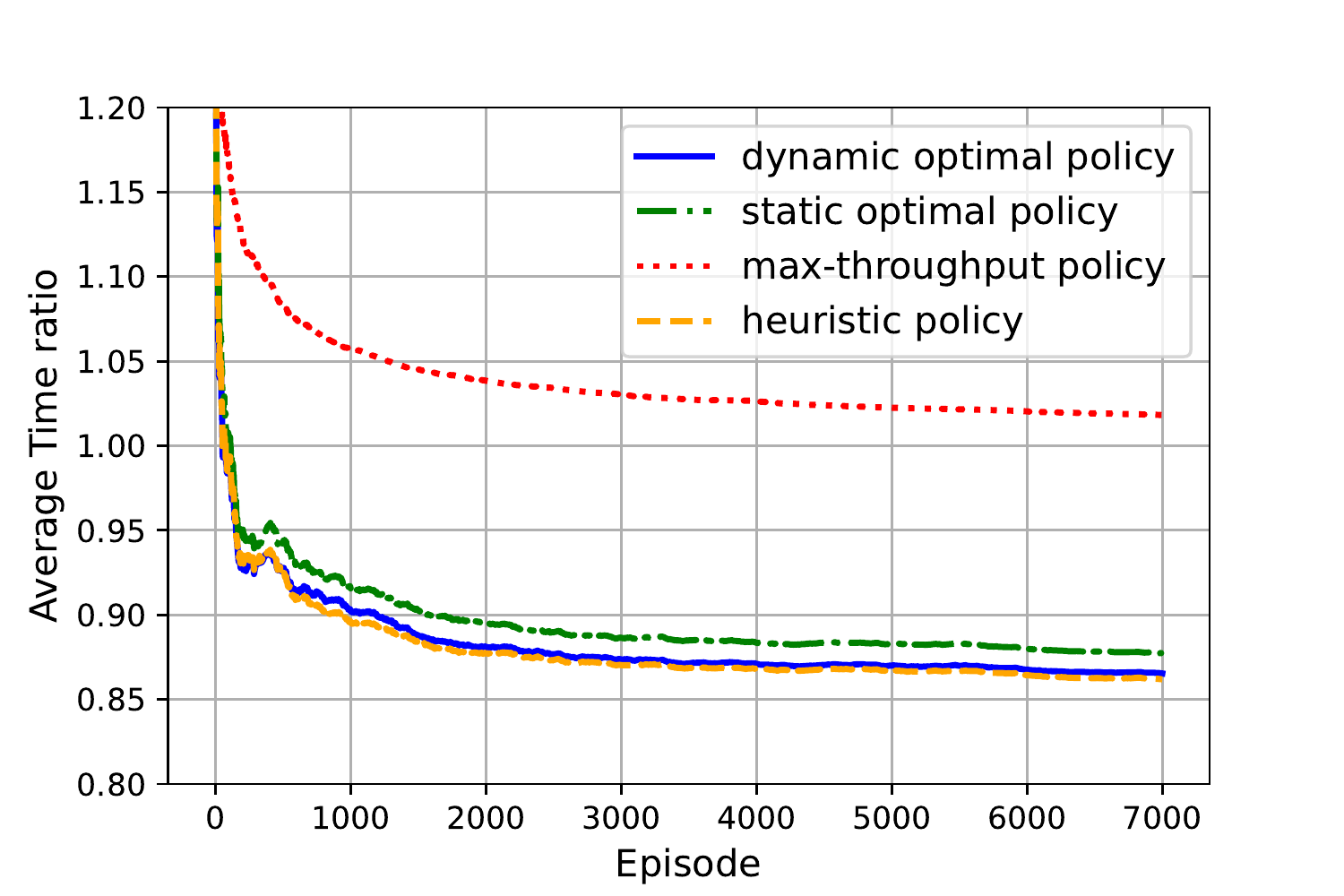}}
     \vspace{-1mm}
     \caption{Effect of switching delay on the performance of (1) dynamic optimal, static optimal, and max-throughput policies in the offline setting with or without switching delay; (2) above three policies and heuristic policy in the online setting with switching delay over the lossy channels.}
     \label{fig:switching_cost_section}\vspace{-4mm}
\end{figure*}
%----------------------------
In the online file transfer problem, since the file size needs not be fixed and larger file sizes naturally take more time to transmit, it makes sense to normalize our performance metric across the range of file sizes and use max-throughput policy as our baseline policy. We define our metrics as \emph{average time ratio} and \emph{average throughput}. For an arbitrary sequence of files $\{F^k\}_{k \in \Z_+}$, the average time ratio at the $K$-th episode is defined as 
\begin{equation}\label{eqn:avg_time_ratio}
    \frac{1}{K}\sum_{k=1}^K{T(\pi^k,F^k)}/{\E[T(i^*,F^k)]},
\end{equation}
and the average throughput is represented as 
\begin{equation}\label{eqn:avg_throughput}
    \frac{1}{K}\sum_{k=1}^K F^k/T(\pi^k,F^k).
\end{equation}
Here, $T(\pi^k,F^k)$ is the measured transfer time of a file of size $F^k$ applying the policy $\pi^k$ at the $k$-th episode. Policy $\pi^k$ is based on the estimated parameter, which is updated by the SU on the fly, as described in Section \ref{subsection:online_implementation}.

In our simulation, we generate $7000$ files from $(0,7]$ (Mb) uniformly at random to be used in Algorithm \ref{alg:online_algorithm}. The simulation is repeated $500$ times to ensure stable results. 
We first observe the bottom row in Figure \ref{fig:online_simulation_gain}. The max-throughput policy achieves the largest average throughput while, counter-intuitively, has the longest transfer time in all channel cases. The reason is that the max-throughput policy computed by the SU is decided by the estimated parameters and can be the inferior policy, resulting in lower average throughput initially, which is an effect of imperfect knowledge of channel parameters. As time goes on, we can see the red curve eventually exceeds all other curves because the SU will eventually learn all the channel parameters well.

Next we focus on the average time ratio in the top row of Figure \ref{fig:online_simulation_gain}. We first observe that all curves eventually flatten out, signifying the convergence of Algorithm \ref{alg:online_algorithm}. In the gradual case, the average time ratio is above $95\%$ for all three policies, implying that they don't obtain much reduction in time and the max-throughput channel is good to access when it is available for most of the time.
However, as shown in the \textit{steep} and \textit{lossy} cases respectively, the dynamic optimal policy and heuristic policy, as well as the static optimal policy, can save over $10\%$ time on average over the baseline. This observation is in line with Corollary \ref{corollary:3.3} and Corollary \ref{corollary:dynamic_lower_upper_bound} since the available probabilities of the max-throughput channel are very small in \textit{steep} and \textit{lossy} cases. Furthermore, the heuristic policy, in addition to keeping the complexity low, achieves similar transfer time to that of the dynamic optimal policy; at the same time performing better than the static optimal policy, as expected from Section \ref{subsection: reusable}.

\subsection{OSA File Transfer with Switching Delay}
In reality, switching from one channel to another also takes some time and could affect the file transfer time if the SU switches too often. In this section, we take the switching delay into consideration. Specifically, we set the switching delay to be $20$ ms as in \cite{nezhad2013semi} (within the $\Delta=100$ ms duration) and simulate the file transfer time in both offline and online settings. 

We first consider the empirical file transfer time over the lossy channels without switching delay in Figure \ref{fig:offline_dynamic_vs_max_throughput}. It indicates that dynamic optimal policy (blue curve) achieves the best performance (the lowest curve) and the simulation results are in line with the theoretical results. The inset shows the average time ratio \eqref{eqn:avg_time_ratio} of each policy (smaller is better), where the average is taken over the file size from $0.1$ Mb to $F$ Mb on the x axis. This is to show the expected file transfer time of each policy compared to the max-throughput policy when the file size falls into a given range that is governed by different applications. Larger $F$ leads to smaller average time ratio for each policy, supporting the discussion after Corollary \ref{corollary:dynamic_lower_upper_bound}. The average time ratio of $F=7$ Mb in Figure \ref{fig:offline_dynamic_vs_max_throughput} is $0.85$ for static optimal policy and $0.83$ for dynamic optimal policy, which is consistent with the top-right plot in Figure \ref{fig:online_simulation_gain}.

On the other hand, switching too frequently penalizes the performance of each policy, and the switching delay may outweigh the time saved by channel switching. The performance gap between dynamic optimal policy and static optimal policy becomes smaller in both offline and online settings when switching delay is taken into account, as shown in Figure \ref{fig:offline_sc} (compared to Figure \ref{fig:offline_dynamic_vs_max_throughput}) and \ref{fig:online_sc} (compared to the top-right plot in Figure \ref{fig:online_simulation_gain}) over the lossy channels. Additionally, our proposed heuristic policy (orange curve) performs slightly better than the dynamic optimal policy and is the best among the four policies in Figure \ref{fig:online_sc}. This depicts that the advantage of our heuristic strategy is not only in terms of computational complexity but also in terms of less channel switching.

\section{Extension to Markovian Channels}\label{File Transfer in Markovian Channels}
The OSA literature has long focused on the throughput-oriented policies for both Bernoulli channels \cite{osa_survey,avner2019multi,zhu2020machine} and Markovian channels, i.e., each channel can be modeled as a two-state discrete-time Markov chain \cite{Zhao_2007,zhao2008opportunistic,liu2010indexability,tekin2011online,liang2019deep}. For the file transfer problem, we have presented static and dynamic policies for Bernoulli channels in Section \ref{static analysis} and \ref{mdp formulation and dynamic policy}. However, for dynamic policies over Markovian channels, the problem is beyond the SSP framework described in Section \ref{mdp formulation and dynamic policy} since the transition probability to the next state (remaining file size) also depends on the past action (last chosen channel), instead of merely the current state and the current action (as in the SSP problem). Augmenting the state space to accommodate for this extra dependency transforms the problem into a POMDP problem, which generally has no known structured solution and is therefore intractable \cite{liu2015online}. For this reason, we extend the file transfer problem to Markovian channels and focus mainly on static policies.
\begin{figure}[!h]
    \centering
    \includegraphics[width=0.8\columnwidth]{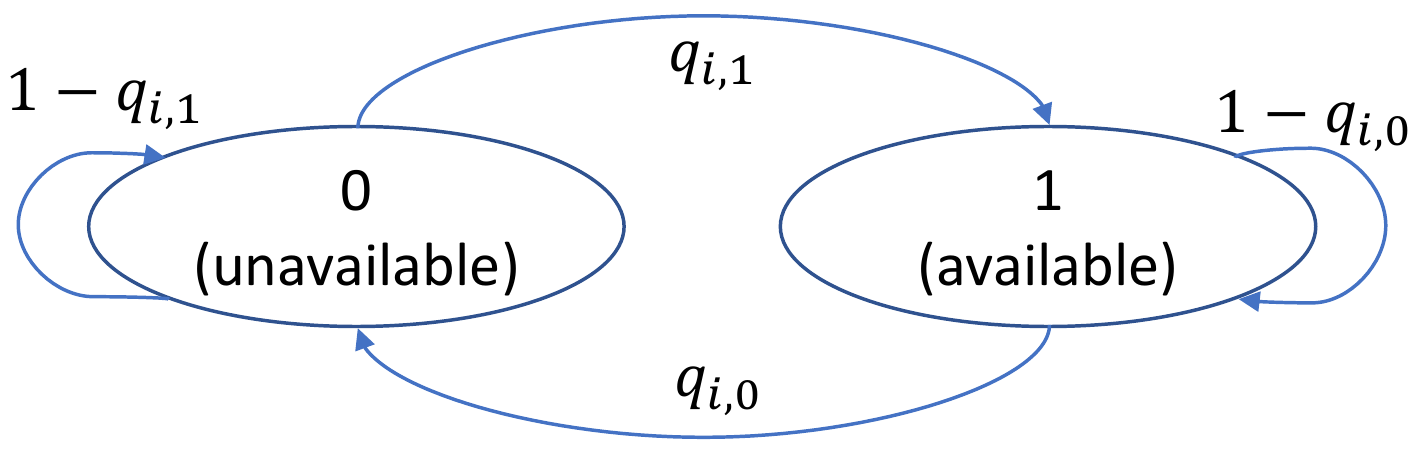}
    \vspace{-2mm}
    \caption{Diagram of the state transition of channel $i\in\cN$ modeled as a two-state Markov chain.}
    \label{fig:diagm_corr}
\end{figure}

In this section, the state of any channel $i\in\cN$ takes the form of a two-state Markov chain $\{Y_i(k)\}_{k\in\N}$. For any time step $k\in\N$, we have $Y_i(k)\in\{0,1\}$, and $P(Y_i(k+1) = 1|Y_i(k)=0) = q_{i,1}$, $P(Y_i(k+1) = 0|Y_i(k)=1) = q_{i,0}$, denoting the transition probabilities, as shown in Figure \ref{fig:diagm_corr}. The stationary distribution of channel $i$ is denoted by $\pi_i \triangleq q_{i,1}/(q_{i,1}+q_{i,0})$. Denote by $c_i\in[0,1]$ an arbitrary probability that channel $i$ is available at the beginning of the file transfer process. Similar to the notations used in Proposition \ref{prop:expected_time_static}, given a file of size $F$, we define $k_i \triangleq \floor{F/\Delta r_i} \in \Z_+$ and $\alpha_i \triangleq  F/ \Delta r_i - k_i \in [0,1)$. Then, we show the closed-form expected file transfer time for static policies in correlated channels as follows.
\begin{proposition}\label{prop:transfer_time_corr}
Given a file of size $F$, the expected transfer time of a static policy in channel $i\in\cN$ is
\begin{equation}\label{eqn:transfer_time_corr}
    \E[T(i,F)] \!=\! \Delta\!\!\left(\!\!\left(k_i+\Char_{\{\alpha_i > 0\}}\!-\!1\right)\frac{q_{i,0}}{q_{i,1}} \!+\! \frac{1\!-\!c_i}{q_{i,1}}\right) \!+\! \frac{F}{r_i}.
\end{equation}
\end{proposition}
\begin{IEEEproof}
Similar to the proof of Proposition \ref{prop:expected_time_static}, we can decompose the file size $F$ as
\begin{equation}\label{eqn:file_size_corr}
    F = k_i \Delta r_i + \alpha_i \Delta r_i
\end{equation}
such that the SU can fully utilize $k_i$ time slots for the file transmission and $\alpha_i$ is the fraction of the $\Delta$ seconds utilized for file transfer toward the end. Initially, the available probability of channel $i$ is $c_i$. The expected waiting time $\E[T_{wait}^1]$ for the first successful transmission of $\Delta r_i$ data is given as
\begin{equation}\label{eqn:waiting_time_1}
\E[T_{wait}^1] = 0\cdot c_i + \sum_{k=1}^{\infty} k(1-c_i)(1-q_{i,1})^{k-1} q_{i,1} = \frac{1-c_i}{q_{i,1}}.
\end{equation}
For the $m$-th transmission of the amount of $\Delta r_i$ data ($m=2,3,\cdots,k_i$), the expected waiting time $\E[T_{wait}^m]$ conditioned on channel $i$ being in state $1$ is
\begin{equation}\label{eqn:waiting_time_2}
\E[T_{wait}^m] = 0 \cdot (1-q_{i,0}) + \sum_{k=1}^{\infty} k q_{i,0}(1-q_{i,1})^{k-1} q_{i,1}  = \frac{q_{i,0}}{q_{i,1}}.
\end{equation}
Note that $\{T_{wait}^m\}_{m=1,2,\cdots,k_i}$ are mutually independent to each other. Let the constant $T_{tran} \triangleq F/r_i = \Delta(k_i+\alpha_i)$ be the total successful transmission time (from \eqref{eqn:file_size_corr}). Then, the transfer time can be written as 
\begin{equation*}
    \begin{split}
        T(i,F) &= \sum_{m=1}^{k_i} \left(T_{wait}^m + \Delta\right) + \Char_{\{\alpha_i > 0\}} T_{wait}^{k_i+1} + \Delta \alpha_i \\
        &= T_{tran} + T_{wait}^1 + \!\!\sum_{m=2}^{k_i+ \Char_{\{\alpha_i > 0\}}}\!\!  T_{wait}^m.
    \end{split}
\end{equation*}
Taking the expectation of the equation above, along with \eqref{eqn:waiting_time_1} and \eqref{eqn:waiting_time_2}, yields \eqref{eqn:transfer_time_corr}.
\end{IEEEproof}
\smallskip The static optimal policy over Markovian channels is then derived from $i_{so}(F) = \argmin_{i\in\cN} \E[T(i,F)]$. By choosing transition probability $q_{i,1} = 1-q_{i,0}$ and initial state $c_i = q_{i,1}$, the Markovian channel reduces to the Bernoulli channel (i.e., $q_{i,1} \equiv p_i$ and $q_{i,0} \equiv 1 - p_i$), and \eqref{eqn:transfer_time_corr} coincides with \eqref{eqn:time_static_policy}.
From \eqref{eqn:transfer_time_corr} we have
\begin{equation}\label{eqn:transfer_time_corr_larger_walds}
\begin{split}
        \E[T(i,F)] &\!=\! \Delta\left((k_i\!+\!\Char_{\{\alpha _i > 0\}})\frac{q_{i,0}}{q_{i,1}} \!+\! \frac{1\!-\!c_i\!-\!q_{i,0}}{q_{i,1}}\right) \!+\! \frac{F}{r_i}\\
        &\!\geq\! \Delta\left((k_i+\alpha_i)\frac{q_{i,0}}{q_{i,1}} + \frac{1-c_i-q_{i,0}}{q_{i,1}}\right) + \frac{F}{r_i}\\
        &\!=\! \frac{F}{r_i\pi_i} + \Delta\frac{1-c_i-q_{i,0}}{q_{i,1}},
\end{split}
\end{equation}
where the inequality comes from $\Char_{\{\alpha_i > 0\}}\geq \alpha_i$. When the initial channel state $c_i$ satisfies the condition $1-c_i-q_i \geq 0$, \eqref{eqn:transfer_time_corr_larger_walds} is lower bounded by $F/r_i\pi_i$, where $r_i\pi_i$ is the standard criterion to choose the max-throughput policy in the single channel \cite{tekin2011online,dai2012efficient}, i.e., $i^* = \argmax_{i\in\cN} r_i\pi_i$. An example for such condition would be that channel $i$ has positive correlation (i.e., $1-q_{i,1}-q_{i,0} > 0$) and is in the stationary regime from the beginning (i.e., $c_i = \pi_i$). Besides, note that when the file size is an integer multiple of $\Delta r_{i^*}$, Proposition \ref{prop:expected_time_static} for Bernoulli channels shows that the expected file transfer time is $F/r_{i^*}p_{i^*}$, while for Markovian channel $i^*$ with positive correlation and $c_{i^*} = \pi_{i^*}$, Proposition \ref{prop:transfer_time_corr} for Markovian channels indicates that $\E[T(i^*,F)] > F/r_{i^*}\pi_{i^*}$ and the equality never holds. Both imply that the max-throughput channel does not necessarily minimize the expected file transfer time.

Now, we analyze the impact of the correlation on the expected file transfer time in each channel. Consider the case where two channels $i$ and $j$ share the same stationary distribution $\pi_i = \pi_j$, initial state probability $c_i=c_j$ and channel rate $r_i = r_j$ but with different correlation, i.e., $q_{i,1} = \beta q_{j,1}$ and $q_{i,0} = \beta q_{j,0}$, where $\beta \in (0,1)$. Then, we know that both channels share the same long-term throughput. However, from Proposition \ref{prop:transfer_time_corr} we have
\begin{equation}
    \E[T(i,F)] - \E[T(j,F)] = \Delta(1-c_i) \left(\frac{1}{q_{i,1}} - \frac{1}{q_{j,1}}\right) > 0,
\end{equation}
which demonstrates that larger correlation (smaller $\beta$) leads to larger expected file transfer time (worse performance).

\section{Conclusion and Future Work}
In this paper, we have developed a theoretical framework for the file transfer problem, where channels are modeled as independent Bernoulli process, to provide the accurate file transfer time for both static and dynamic policies. We pointed out that the max-throughput channel does not always minimize the file transfer time. We demonstrated in our our analysis that our proposed policies can obtain significant reduction in file transfer time over the max-throughput policy for small file sizes or when the max-throughput channel has very high rate but with low available probability, as typically the case in reality. In addition, we have extended the theoretical analysis to Markovian channels and static polices, showing that greater correlation can compromise the performance.

When the wireless devices work in the outside network, the effect of a propagation environment on a radio signal needs to be considered, i.e., Rayleigh fading, Rician fading and Nakagami fading, which leads to the varying channel environment. Our future work includes the extension to the channels with multiple rates such that the rate will be treated as a general random variable sampled from some probability distribution or finite-state Markov chain, instead of Bernoulli random variable or two-state Markov chain studied in this paper.

%% file: main.bbl
\begin{thebibliography}{10}
\providecommand{\url}[1]{#1}
\csname url@samestyle\endcsname
\providecommand{\newblock}{\relax}
\providecommand{\bibinfo}[2]{#2}
\providecommand{\BIBentrySTDinterwordspacing}{\spaceskip=0pt\relax}
\providecommand{\BIBentryALTinterwordstretchfactor}{4}
\providecommand{\BIBentryALTinterwordspacing}{\spaceskip=\fontdimen2\font plus
\BIBentryALTinterwordstretchfactor\fontdimen3\font minus
  \fontdimen4\font\relax}
\providecommand{\BIBforeignlanguage}[2]{{%
\expandafter\ifx\csname l@#1\endcsname\relax
\typeout{** WARNING: IEEEtran.bst: No hyphenation pattern has been}%
\typeout{** loaded for the language `#1'. Using the pattern for}%
\typeout{** the default language instead.}%
\else
\language=\csname l@#1\endcsname
\fi
#2}}
\providecommand{\BIBdecl}{\relax}
\BIBdecl

\bibitem{hu2021opportunistic}
J.~Hu, V.~Doshi, and D.~Y. Eun, ``Opportunistic spectrum access: Does
  maximizing throughput minimize file transfer time?'' in \emph{19th
  International Symposium on Modeling and Optimization in Mobile, Ad hoc, and
  Wireless Networks (WiOpt)}, Philadelphia, PA, USA, 2021.

\bibitem{osa_survey}
K.~Zaheer, M.~Othman, M.~H. Rehmani, and T.~Perumal, ``A survey of
  decision-theoretic models for cognitive internet of things (ciot),''
  \emph{IEEE Access}, vol.~6, pp. 22\,489--22\,512, 2018.

\bibitem{federal_2020}
\BIBentryALTinterwordspacing
``{FCC} increases unlicensed wireless operations in tv white spaces,'' Dec
  2020. [Online]. Available:
  \url{https://www.fcc.gov/document/fcc-increases-unlicensed-wireless-operations-tv-white-spaces-0}
\BIBentrySTDinterwordspacing

\bibitem{802.22_standard}
``{IEEE} standard for information technology,'' \emph{IEEE Std 802.22-2019},
  pp. 1--1465, 2020.

\bibitem{naparstek2018deep}
O.~Naparstek and K.~Cohen, ``Deep multi-user reinforcement learning for
  distributed dynamic spectrum access,'' \emph{IEEE Transactions on Wireless
  Communications}, vol.~18, no.~1, pp. 310--323, 2018.

\bibitem{avner2019multi}
O.~Avner and S.~Mannor, ``Multi-user communication networks: A coordinated
  multi-armed bandit approach,'' \emph{IEEE/ACM Transactions on Networking},
  vol.~27, no.~6, pp. 2192--2207, 2019.

\bibitem{nasir2019multi}
Y.~S. Nasir and D.~Guo, ``Multi-agent deep reinforcement learning for dynamic
  power allocation in wireless networks,'' \emph{IEEE IEEE Journal on Selected
  Areas in Communications.}, vol.~37, no.~10, pp. 2239--2250, 2019.

\bibitem{liu2010indexability}
K.~Liu and Q.~Zhao, ``Indexability of restless bandit problems and optimality
  of whittle index for dynamic multichannel access,'' \emph{IEEE Transactions
  on Information Theory}, vol.~56, no.~11, pp. 5547--5567, 2010.

\bibitem{yadav2022deep}
M.~A. Yadav, Y.~Li, G.~Fang, and B.~Shen, ``Deep q-network based reinforcement
  learning for distributed dynamic spectrum access,'' in \emph{2022 IEEE 2nd
  International Conference on Computer Communication and Artificial
  Intelligence (CCAI)}, 2022, pp. 1--6.

\bibitem{Online_Sequential_2}
P.~Yang, B.~Li, J.~Wang, X.~Li, Z.~Du, Y.~Yan, and Y.~Xiong, ``Online
  sequential channel accessing control: A double exploration vs. exploitation
  problem,'' \emph{IEEE Transactions on Wireless Communications}, vol.~14,
  no.~8, pp. 4654--4666, 2015.

\bibitem{zuo2019obp}
J.~Zuo, X.~Zhang, and C.~Joe-Wong, ``Observe before play: Multi-armed bandit
  with pre-observations,'' in \emph{Proceedings of the AAAI Conference on
  Artificial Intelligence}, New York, USA, 2020.

\bibitem{monemian2015optimum}
M.~Monemian, M.~Mahdavi, and M.~J. Omidi, ``Optimum sensor selection based on
  energy constraints in cooperative spectrum sensing for cognitive radio sensor
  networks,'' \emph{IEEE Sensors Journal}, vol.~16, no.~6, pp. 1829--1841,
  2015.

\bibitem{gan2018cost}
C.~Gan, R.~Zhou, J.~Yang, and C.~Shen, ``Cost-aware learning and optimization
  for opportunistic spectrum access,'' \emph{IEEE Transactions on Cognitive
  Communications and Networking}, vol.~5, no.~1, pp. 15--27, 2018.

\bibitem{liang2019deep}
L.~Liang, H.~Ye, G.~Yu, and G.~Y. Li, ``Deep-learning-based wireless resource
  allocation with application to vehicular networks,'' \emph{Proceedings of the
  IEEE}, vol. 108, no.~2, pp. 341--356, 2019.

\bibitem{Low_Latency_Towards_5G}
I.~Parvez, A.~Rahmati, I.~Guvenc, A.~I. Sarwat, and H.~Dai, ``A survey on low
  latency towards 5g: Ran, core network and caching solutions,'' \emph{IEEE
  Communications Surveys \& Tutorials}, vol.~20, no.~4, pp. 3098--3130, 2018.

\bibitem{hussain2018autonomous}
R.~Hussain and S.~Zeadally, ``Autonomous cars: Research results, issues, and
  future challenges,'' \emph{IEEE Communications Surveys \& Tutorials},
  vol.~21, no.~2, pp. 1275--1313, 2018.

\bibitem{paul2019spectrum}
A.~Paul, P.~Kunarapu, A.~Banerjee, and S.~P. Maity, ``Spectrum sensing in
  cognitive vehicular networks for uniform mobility model,'' \emph{IET
  Communications}, vol.~13, no.~19, pp. 3127--3134, 2019.

\bibitem{8024608}
G.~V. Rossi and K.~K. Leung, ``Optimised csma/ca protocol for safety messages
  in vehicular ad-hoc networks,'' in \emph{2017 IEEE Symposium on Computers and
  Communications (ISCC)}, Heraklion, Greece.

\bibitem{sodagari2018technologies}
S.~Sodagari, B.~Bozorgchami, and H.~Aghvami, ``Technologies and challenges for
  cognitive radio enabled medical wireless body area networks,'' \emph{IEEE
  Access}, vol.~6, pp. 29\,567--29\,586, 2018.

\bibitem{Srikant_2018}
H.~Gupta, A.~Eryilmaz, and R.~Srikant, ``Low-complexity, low-regret link rate
  selection in rapidly-varying wireless channels,'' in \emph{IEEE INFOCOM},
  Honolulu, HI, USA, 2018.

\bibitem{Srikant_2019}
------, ``Link rate selection using constrained thompson sampling,'' in
  \emph{IEEE INFOCOM}, Paris, France, 2019.

\bibitem{zhu2020machine}
P.~Zhu, J.~Li, D.~Wang, and X.~You, ``Machine-learning-based opportunistic
  spectrum access in cognitive radio networks,'' \emph{IEEE Wireless
  Communications}, vol.~27, no.~1, pp. 38--44, 2020.

\bibitem{almasri2021managing}
M.~Almasri, A.~Mansour, C.~Moy, A.~Assoum, D.~Le~Jeune, and C.~Osswald,
  ``Managing single or multi-users channel allocation for the priority
  cognitive access,'' in \emph{2020 28th European Signal Processing Conference
  (EUSIPCO)}, Amsterdam, NL.

\bibitem{wu2014learning}
Y.~Wu, F.~Hu, S.~Kumar, Y.~Zhu, A.~Talari, N.~Rahnavard, and J.~D. Matyjas, ``A
  learning-based qoe-driven spectrum handoff scheme for multimedia
  transmissions over cognitive radio networks,'' \emph{IEEE Journal on Selected
  Areas in Communications}, vol.~32, no.~11, pp. 2134--2148, 2014.

\bibitem{cao2017dynamic}
H.~Cao, H.~Tian, J.~Cai, A.~S. Alfa, and S.~Huang, ``Dynamic load-balancing
  spectrum decision for heterogeneous services provisioning in multi-channel
  cognitive radio networks,'' \emph{IEEE Transactions on Wireless
  Communications}, vol.~16, no.~9, pp. 5911--5924, 2017.

\bibitem{dimitriou2018stable}
I.~Dimitriou and N.~Pappas, ``Stable throughput and delay analysis of a random
  access network with queue-aware transmission,'' \emph{IEEE Transactions on
  Wireless Communications}, vol.~17, no.~5, pp. 3170--3184, 2018.

\bibitem{huang2019dynamic}
X.-L. Huang, X.-W. Tang, and F.~Hu, ``Dynamic spectrum access for multimedia
  transmission over multi-user, multi-channel cognitive radio networks,''
  \emph{IEEE Transactions on Multimedia}, vol.~22, no.~1, pp. 201--214, 2019.

\bibitem{iqbal2021enhanced}
A.~Iqbal, R.~Hussain, A.~Shakeel, I.~L. Khan, M.~A. Javed, Q.~U. Hasan, B.~M.
  Lee, and S.~A. Malik, ``Enhanced spectrum access for qos provisioning in
  multi-class cognitive d2d communication system,'' \emph{IEEE Access}, vol.~9,
  pp. 33\,608--33\,624, 2021.

\bibitem{tschabitscher_2020}
\BIBentryALTinterwordspacing
Lifewire, ``Ever wonder what makes email files so large?'' 2020. [Online].
  Available:
  \url{https://www.lifewire.com/what-is-the-average-size-of-an-email-message-1171208}
\BIBentrySTDinterwordspacing

\bibitem{Zhao_2007}
Q.~Zhao, L.~Tong, A.~Swami, and Y.~Chen, ``Decentralized cognitive mac for
  opportunistic spectrum access in ad hoc networks: A pomdp framework,''
  \emph{IEEE Journal on selected areas in communications}, vol.~25, no.~3, pp.
  589--600, 2007.

\bibitem{liu2015online}
Y.~Liu and M.~Liu, ``An online approach to dynamic channel access and
  transmission scheduling,'' in \emph{Proceedings of the 16th ACM International
  Symposium on Mobile Ad Hoc Networking and Computing}, Hangzhou, China, 2015.

\bibitem{zhao2008opportunistic}
Q.~Zhao, S.~Geirhofer, L.~Tong, and B.~M. Sadler, ``Opportunistic spectrum
  access via periodic channel sensing,'' \emph{IEEE Transactions on Signal
  Processing}, vol.~56, no.~2, pp. 785--796, 2008.

\bibitem{tekin2011online}
C.~Tekin and M.~Liu, ``Online learning in opportunistic spectrum access: A
  restless bandit approach,'' in \emph{IEEE INFOCOM}, Shanghai, China, 2011.

\bibitem{dai2012efficient}
W.~Dai, Y.~Gai, and B.~Krishnamachari, ``Efficient online learning for
  opportunistic spectrum access,'' in \emph{IEEE INFOCOM}, Orlando, FL, USA,
  2012.

\bibitem{NEURIPS2019_2edfeadf}
Y.~H. Jung and A.~Tewari, ``Regret bounds for thompson sampling in episodic
  restless bandit problems,'' in \emph{Advances in Neural Information
  Processing Systems}, Vancouver, Canada, 2019.

\bibitem{wang2018deep}
S.~Wang, H.~Liu, P.~H. Gomes, and B.~Krishnamachari, ``Deep reinforcement
  learning for dynamic multichannel access in wireless networks,'' \emph{IEEE
  Transactions on Cognitive Communications and Networking}, vol.~4, no.~2, pp.
  257--265, 2018.

\bibitem{mohamedou2017bayesian}
A.~Mohamedou, A.~Sali, B.~Ali, M.~Othman, and H.~Mohamad, ``Bayesian inference
  and fuzzy inference for spectrum sensing order in cognitive radio networks,''
  \emph{Transactions on Emerging Telecommunications Technologies}, vol.~28,
  no.~1, p. e2916, 2017.

\bibitem{talebi2017stochastic}
M.~S. Talebi, Z.~Zou, R.~Combes, A.~Proutiere, and M.~Johansson, ``Stochastic
  online shortest path routing: The value of feedback,'' \emph{IEEE
  Transactions on Automatic Control}, vol.~63, no.~4, pp. 915--930, 2017.

\bibitem{tan2022cooperative}
X.~Tan, L.~Zhou, H.~Wang, Y.~Sun, H.~Zhao, B.-C. Seet, J.~Wei, and V.~C. Leung,
  ``Cooperative multi-agent reinforcement learning based distributed dynamic
  spectrum access in cognitive radio networks,'' \emph{IEEE Internet of Things
  Journal}, 2022.

\bibitem{pei2009much}
Y.~Pei, Y.-C. Liang, K.~C. Teh, and K.~H. Li, ``How much time is needed for
  wideband spectrum sensing?'' \emph{IEEE Transactions on Wireless
  Communications}, vol.~8, no.~11, pp. 5466--5471, 2009.

\bibitem{toma2020estimation}
O.~H. Toma, M.~Lopez-Benitez, D.~K. Patel, and K.~Umebayashi, ``Estimation of
  primary channel activity statistics in cognitive radio based on imperfect
  spectrum sensing,'' \emph{IEEE Transactions on Communications}, vol.~68,
  no.~4, pp. 2016--2031, 2020.

\bibitem{bicket2005bit}
J.~C. Bicket, ``Bit-rate selection in wireless networks,'' Master's thesis,
  Massachusetts Institute of Technology, 2005.

\bibitem{bertsekas2019reinforcement}
D.~Bertsekas, \emph{Reinforcement Learning and Optimal Control}.\hskip 1em plus
  0.5em minus 0.4em\relax Athena Scientific, 2019.

\bibitem{della2017exact}
F.~Della~Croce, F.~Salassa, and R.~Scatamacchia, ``An exact approach for the
  0--1 knapsack problem with setups,'' \emph{Computers \& Operations Research},
  vol.~80, pp. 61--67, 2017.

\bibitem{conforti2014integer}
M.~Conforti, G.~Cornu{\'e}jols, G.~Zambelli \emph{et~al.}, \emph{Integer
  programming}.\hskip 1em plus 0.5em minus 0.4em\relax Springer, 2014, vol.
  271.

\bibitem{bertsekas2012dynamic}
D.~Bertsekas, \emph{Dynamic programming and optimal control: Volume I}.\hskip
  1em plus 0.5em minus 0.4em\relax Athena scientific, 2012, vol.~1.

\bibitem{mantuano_2016}
\BIBentryALTinterwordspacing
A.~Mantuano, ``File size basics,'' 2016. [Online]. Available:
  \url{https://techdocs.blogs.brynmawr.edu/5523}
\BIBentrySTDinterwordspacing

\bibitem{gai2012combinatorial}
Y.~Gai, B.~Krishnamachari, and R.~Jain, ``Combinatorial network optimization
  with unknown variables: Multi-armed bandits with linear rewards and
  individual observations,'' \emph{IEEE/ACM Transactions on Networking},
  vol.~20, no.~5, pp. 1466--1478, 2012.

\bibitem{chen2013combinatorial}
W.~Chen, Y.~Wang, and Y.~Yuan, ``Combinatorial multi-armed bandit: General
  framework and applications,'' in \emph{International conference on machine
  learning}, Atlanta, GA, USA, 2013.

\bibitem{MaherMiltenbergerPedrosoRehfeldtSchwarzSerrano2016}
S.~Maher, M.~Miltenberger, J.~P. Pedroso, D.~Rehfeldt, R.~Schwarz, and
  F.~Serrano, ``{PySCIPOpt}: Mathematical programming in python with the {SCIP}
  optimization suite,'' in \emph{Mathematical Software {\textendash} {ICMS}
  2016}.\hskip 1em plus 0.5em minus 0.4em\relax Springer International
  Publishing, 2016, pp. 301--307.

\bibitem{nezhad2013semi}
M.~A. Nezhad, L.~Cerd{\`a}-Alabern, B.~Bellalta, and M.~G. Zapata, ``A
  semi--dynamic, game based and interference aware channel assignment for
  multi--radio multi--channel wireless mesh networks,'' \emph{International
  Journal of Ad Hoc and Ubiquitous Computing}, vol.~14, no.~3, pp. 200--213,
  2013.

\end{thebibliography}
